\documentclass{article}

\usepackage{amsmath}
\usepackage{amsthm}
\usepackage[psamsfonts]{amsfonts}
\usepackage{amssymb}
\usepackage[all]{xy}

\DeclareMathSymbol{\rightrightarrows}  {\mathrel}{AMSa}{"13}

\def\Ob{\operatorname{Ob}}
\def\Mor{\operatorname{Mor}}

\def\Pre{\operatorname{Pre}}
\def\Ho{\operatorname{Ho}}
\def\Ex{\operatorname{Ex}}
\def\pb{\operatorname{pb}}

\catcode`\@=11\def\varholim@#1#2{\mathop{\vtop{\ialign{##\crcr
 \hfil$#1\m@th\operator@font holim$\hfil\crcr
 \noalign{\nointerlineskip\kern\ex@}#2#1\crcr
 \noalign{\nointerlineskip\kern-\ex@}\crcr}}}}
\def\hocolim{\mathpalette\varholim@\rightarrowfill@} 
\def\hoinvlim{\mathpalette\varholim@\leftarrowfill@}
\catcode`\@=\active

\newtheorem{theorem}{Theorem}
\newtheorem{lemma}[theorem]{Lemma}
\newtheorem{corollary}[theorem]{Corollary}

\theoremstyle{definition}

\newtheorem{example}[theorem]{Example}

\begin{document}

\title{Fibred sites and stack cohomology}
 
\author{J.F. Jardine}

 
\date{June 21, 2004}

\maketitle

\section*{Introduction}

A stack $G$ is traditionally defined to be a pseudofunctor on
a Grothendieck site $\mathcal{C}$ which takes values in groupoids,
and which satisfies the effective descent condition. The effective descent
condition specifies that the objects of $G$ satisfy a
pseudo-functorial sheaf condition. With this description in hand, one
can form the Grothendieck construction, here denoted by $\mathcal{C}
/ G$, for the stack $G$ and let it inherit a topology from the
ambient site $\mathcal{C}$, so that $\mathcal{C} / G$
acquires the structure of a Grothendieck site. Then one can speak of
sheaves on this site, and stack cohomology of $G$ with coefficients in
a sheaf $F$ on $\mathcal{C} / G$ is the cohomology of the
site with coefficients in $F$ in the standard way.

That said, the connection between this definition of stack cohomology
and the geometry of the stack $G$ is a bit tenuous, at least
apparently, and it has historically been rather awkward to relate this
invariant to other standard sheaf-theoretic invariants.

The general approach to stacks (and higher stacks) has changed a great
deal in recent years, because we now understand that they are homotopy
theoretic objects. A stack $G$ is now thought of, most generally, as a
presheaf of groupoids on a site $\mathcal{C}$ which is fibrant with
respect to a nicely defined model structure on the category of
presheaves of groupoids on $\mathcal{C}$

More explicitly, one says that a functor $G \to H$ between presheaves
of groupoids is a weak equivalence (respectively fibration) if the
induced map $BG \to BH$ is a local weak equivalence (respectively
global fibration) in the standard model structure on the
category of simplicial presheaves on $\mathcal{C}$. Thus,
$G$ is a stack if and only if $BG$ is a globally fibrant simplicial
presheaf. This description of stacks was a major conceptual
breakthrough which was initiated by Joyal and Tierney \cite{JT0} in
the case of sheaves of groupoids and was completed by Hollander
\cite{H} for presheaves of groupoids. Stack completion becomes a
fibrant model in this setup, and it is now well understood that path
components (or isomorphism classes) in the global sections of a stack
$G$ are in bijective correspondence with the set $[\ast,BG]$ of
morphisms in the homotopy category of simplicial presheaves. This
gives a rather striking generalization of the early result that
identified the homotopy invariants $[\ast,BH]$ arising from sheaves of
groups $H$ with isomorphism classes of $H$-torsors \cite{UH}. We also now
understand what the higher order analogues of $H$-torsors should be,
and a homotopy theoretic (and geometric) identification of these
higher order torsors has been achieved \cite{JL}.

This paper brings stack cohomology into this arena, by giving an
homotopy theoretic description of the invariant in terms of presheaves
of groupoids. One of the more important consequences of this approach
is that one can then show that the new cohomology theory for
presheaves of groupoids is homotopy invariant.

In fact, one generalizes the traditional description of the site
$\mathcal{C} / G$ fibred over a stack $G$ even further, to that of the
site $\mathcal{C} / A$ fibred over a presheaf of categories $A$.  This
seems like a strange thing to do at first, but the concept is painless
to both define and manipulate.  This expanded notion specializes to
fibred site constructions that are in standard use, including the
usual sites fibred over diagrams of schemes, and hence over simplicial
schemes in standard geometric settings. It is also interesting to observe
that the idea has non-trivial content even in the case where $A$
consists only of a presheaf of objects.

Simplicial presheaves $X$ for the site $\mathcal{C} / A$ take
the form of enriched contravariant diagrams defined on $A$ and taking
values in simplicial sets, and as such naturally determine homotopy
colimits
\begin{equation*}
\hocolim_{A^{op}}\ X \to BA^{op}
\end{equation*}
over the nerve $BA^{op}$ of the opposite category $A^{op}$, or
equivalently over $BA$.  The homotopy theory of simplicial presheaves
on the fibred site $\mathcal{C} / A$ is actually a type of coarse
equivariant theory of $A^{op}$-diagrams --- one says ``coarse'' because
this is an enriched version of the old Bousfield-Kan theory for
diagrams of simplicial sets \cite{BK}.

In the case when $A$ is a presheaf of groupoids $G$, this assignment
of homotopy colimits determines an equivalence of homotopy categories
\begin{equation}\label{neweq 1}
\Ho(s\Pre(\mathcal{C} / G)) \simeq \Ho(s\Pre(\mathcal{C})
/ BG^{op})
\end{equation}
which generalizes the known relationship \cite{GJ} between diagrams of
simplicial sets defined on a groupoid $H$ and that of simplicial sets
over $BH$. This identification gives the homotopy invariance, because
the homotopy category of simplicial presheaves over $BG^{op}$ is
insensitive to the homotopy type of the presheaf of groupoids $G$ up
to equivalence.

It is a consequence of the equivalence (\ref{neweq 1}) that any
functor $G \to H$ of presheaves of groupoids which is a local weak
equivalence induces an adjoint equivalence of homotopy categories
\begin{equation}\label{neweq 2}
\Ho(s\Pre(\mathcal{C}/G)) \simeq
\Ho(s\Pre(\mathcal{C}/H))
\end{equation}
With a little care (so that you don't spend a long time doing it),
this adjoint equivalence can be parlayed into an adjoint equivalence
\begin{equation}\label{neweq 2.5}
\Ho(s_{\ast}\Pre(\mathcal{C}/G)) \simeq
\Ho(s_{\ast}\Pre(\mathcal{C}/H))
\end{equation}
for pointed simplicial presheaves, and then to adjoint equivalences of
stable homotopy categories
\begin{equation}\label{neweq 3}
\Ho(\mathbf{Spt}(\mathcal{C} / G)) \simeq
\Ho(\mathbf{Spt}(\mathcal{C} / H))
\end{equation}
and
\begin{equation}\label{neweq 4}
\Ho(\mathbf{Spt}_{\Sigma}(\mathcal{C} / G)) \simeq
\Ho(\mathbf{Spt}_{\Sigma}(\mathcal{C} / H))
\end{equation}
for presheaves of spectra and presheaves of symmetric spectra,
respectively. 

Note the level of generality: these results hold over
arbitrary small Grothen\-dieck sites $\mathcal{C}$.
One fully expects that this pattern can be replicated for categories
of module spectra and for various derived categories of chain
complexes, as the need arises. The development given in this paper
ends with the symmetric spectrum result.

The overall aim of this paper is to introduce a very general new
construction, namely the site fibred over a presheaf of categories,
and to give some of its applications for presheaves of groupoids.  The
statements which are listed as theorems are Theorem \ref{thm 30},
which establishes the equivalence (\ref{neweq 1}), and Theorem
\ref{thm 41}, which gives (\ref{neweq 3}). The equivalence (\ref{neweq
2}) appears here as Corollary \ref{cor 32}, and (\ref{neweq 2.5}) is
Corollary \ref{cor 37}. The equivalence (\ref{neweq 4}) is a
consequence of Theorem \ref{thm 41}, and is formally stated in
Corollary \ref{cor 45}.

My personal impression is that the homotopy invariance statements will
turn out to be quite important in applications, as one now has the
ability to define stack cohomology via a construction which comes
directly out of a representing presheaf of groupoids without passing
to any form of either associated sheaf or stack completion. The first
example that comes to mind for which this may be of some use is in the
applications of the cohomology of the presheaf of formal group laws on
the flat site.
\medskip

I would like to thank the American Institute of Mathematics for its
hospitality and support during the week of the workshop ``Theory of
motives, homotopy theory of varieties and dessins d'enfants'', held at
AIM April 23--26, 2004. The appearance of this paper is a direct
result of my participation in that conference.

\tableofcontents
\vfill\eject

\section{Fibred sites}

Suppose that $\mathcal{C}$ is a small Grothendieck site, and that $A$
is a presheaf of categories on $\mathcal{C}$. 

The category $\mathcal{C} / A$ has objects consisting of all
pairs $(U,x)$ where $U$ is an object of $\mathcal{C}$ and $x$ is an
element of the set of section $\Ob(A)(U)$ of the presheaf $\Ob(A)$ of
objects of $A$. We can and will alternatively think of $x$ as a
presheaf morphisms $x: U \rightarrow \Ob(A)$. A morphism
$(\alpha,f): (V,y) \rightarrow (U,x)$ in this category is a pair
consisting of a morphism $\alpha: V \to U$ of $\mathcal{C}$ together
with a morphism $f: y \rightarrow \alpha^{\ast}(x)$ of
$A(U)$.

Given another morphism
$(\gamma,g): (W,z) \to (V,y)$, 
the composite
$(\alpha,f)(\gamma,g)$ is defined by
\begin{equation*}
(\alpha,f)(\gamma,g) =
(\alpha\gamma,\gamma^{\ast}(f)g),
\end{equation*}
where the composite
\begin{equation*}
z \xrightarrow{g} \gamma^{\ast}(y)
\xrightarrow{\gamma^{\ast}(f)}
\gamma^{\ast}\alpha^{\ast}(x) = (\alpha\gamma)^{\ast}(x)
\end{equation*}
is defined by the usual sort of convolution --- this category is the
result of applying the 
Grothendieck construction to the diagram of categories represented by $A$.

There is a canonical forgetful functor $\pi: \mathcal{C} / A
\rightarrow \mathcal{C}$ which is defined by sending the object
$(U,x)$ to the object $U$ of $\mathcal{C}$. Observe that any sieve $R
\subset \hom(\ ,(U,x))$ of $\mathcal{C} / A$ is mapped to a sieve
$\pi(R) \subset \hom(\ ,U)$ under the functor $\pi$.

If $S$ is a sieve for $U \in \mathcal{C}$ write $\pi^{-1}(S)$ for the
collection of all morphisms $(\alpha,f)$ with $\alpha \in
S$. The covering sieves of $\mathcal{C} / A$ are the sieves
of the form $\pi^{-1}S$ for covering sieves $S$ of $\mathcal{C}$. If
$R$ contains a covering sieve $\pi^{-1}S$, then $R = \pi^{-1}(S')$ for
some covering sieve of $S$. Note that $\pi^{-1}S$ is the smallest sieve
containing all morphisms $(\alpha,1)$ with $\alpha \in S$.

\begin{lemma}
The collection of covering sieves $R$ for $\mathcal{C} / A$
satisfy the axioms for a Grothendieck topology.
\end{lemma}

\begin{proof}
There is a relation $(\alpha,f)^{-1}\pi^{-1}S =
\pi^{-1}\alpha^{-1}S$, so the covering sieves of $\mathcal{C}
/ Y$ are closed under pullback.

Suppose that $S$ is a covering sieve for $U$ and that $S_{\alpha}$ is
a choice of covering sieve for $V$ for each $\alpha: V \to U$ in
$S$. Then the local character axiom for the site $\mathcal{C}$ implies
that the collection of all maps $W' \to U$ which factor through a
composite
\begin{equation*}
W \xrightarrow{\beta} V \xrightarrow{\alpha} U
\end{equation*}
with $\alpha \in S$ and $\beta \in S_{\alpha}$ is a covering sieve for $U$.

Suppose that $R,R'$ are sieves for $(U,x)$ and that $R$ is
covering. Suppose further that $(\alpha,f)^{-1}(R')$ is covering
for all $(\alpha,f) \in R$. Suppose that $R = \pi^{-1}S$. Then
$(\alpha,1)^{-1}(R')$ is a covering sieve for each $\alpha: V \to U$
in $S$, and so there is a covering sieve $S_{\alpha}$ for $U$ such
that $\alpha^{-1}R'$ contains all morphisms $(\gamma,1)$
with $\gamma \in S_{\alpha}$. It follows that the collection of all
morphisms of the form $(\zeta,1)$ in $R'$
defines a covering sieve of $\mathcal{C}$ for $U$, so that
$R'$ is covering.

The (trivial) sieve of all morphisms $(\alpha,f)$
is $\pi^{-1}S$, where $S$ is the sieve of all morphisms $V \to U$ in
$\mathcal{C}$. All trivial sieves are therefore covering.
\end{proof}

The site $\mathcal{C} / A$ will be called the {\it fibred
site} for the presheaf of categories $A$.  Recall that every
Grothendieck site $\mathcal{D}$ determines a model structure on its
associated category $s\Pre(\mathcal{D})$ of simplicial presheaves, for
which the cofibrations are the monomorphisms, the weak equivalences
are the local weak equivalences, and the fibrations are the global
fibrations. One of the primary goals of this paper is to analyze the
corresponding model structures on $s\Pre(\mathcal{C} / A)$
for the fibred sites in important special cases. Sites fibred over
presheaves of groupoids will be of fundamental interest in the later
sections of this paper. This construction is quite general, and
specializes to many other well known standard examples, as the
following preliminary list of examples is meant to demonstrate.

\begin{example}\label{ex 2}
Suppose that $I$ is a small category and that $Y: I \to
\Pre(\mathcal{C})$ is an $I$-diagram in the category of presheaves on
$\mathcal{C}$ defined by $i \mapsto Y_{i}$. Write
$\mathcal{C} / Y$ for the category whose objects consist of
presheaf morphisms (ie. sections) $x: U \to Y_{i}$, and whose morphisms
are commutative diagrams of presheaf morphisms
\begin{equation*}
\xymatrix{
V \ar[r]^{\alpha} \ar[d]_{y} & U \ar[d]^{x} \\
Y_{i} \ar[r]_{\theta_{\ast}} & Y_{j}
}
\end{equation*}
where $\theta: i \to j$ is a morphism of $I$. Denote such a morphism
by $(\alpha,\theta)$.

There is a presheaf of categories $EY$ which is defined
by setting $EY(U)$ to be the translation category for the functor
$Y(U): I \to \mathbf{Set}$, where $Y(U)_{i} = Y_{i}(U)$.The objects of
$EY(U)$ consist of all pairs $(i,x)$ where $x$ is an element of
$Y_{i}(U)$. Equivalently, $x$ is a presheaf map $U \to Y_{i}$. If
$y: V \to Y_{j}$ is an object of $EY(V)$ and $\alpha: V \to U$ is a
morphism of $\mathcal{C}$, then $\alpha^{\ast}(x)$ is the composite
\begin{equation*}
V \xrightarrow{\alpha} U \xrightarrow{x} Y_{j},
\end{equation*}
so that a morphism $f: (i,y) \to (j,\alpha^{\ast}(x))$ in
$EY(V)$ is morphism $\theta: i \to j$ such that the diagram
\begin{equation*}
\xymatrix{
V \ar[r]^{\alpha} \ar[d]_{y} & U \ar[d]^{x} \\
Y_{i} \ar[r]_{\theta_{\ast}} & Y_{j}
}
\end{equation*}
commutes.  The category $\mathcal{C} / Y$, as defined above,
therefore coincides with the category $C / EY$.
\end{example}

\begin{example}\label{ex 3}
The previous example specializes to the ``standard'' description of
the site $\mathcal{C} / X$ for a simplicial presheaf $X$. 

Suppose that $X$ is a simplicial presheaf and that $Z$ is a globally
fibrant simplicial presheaf on the site $\mathcal{C} / X$.
For each $n$, there is a presheaf $1_{X_{n}}$ which is represented by
the identity on $X_{n}$ in the sense that the sections corresponding
to $y: V \to X_{k}$ are the commutative diagrams of
$\mathcal{C}$-presheaf maps
\begin{equation}\label{eq 1}
\xymatrix{
V \ar[r] \ar[d]_{y} & X_{n} \ar[d]^{1} \\
X_{k} \ar[r]_{\theta^{\ast}} & X_{n}
}
\end{equation}
where $\theta^{\ast}$ is a simplicial structure map. The maps
\begin{equation*}
\xymatrix{
X_{k} \ar[r]^{\theta^{\ast}} \ar[d]_{1} & X_{n} \ar[d]^{1} \\
X_{k} \ar[r]_{\theta^{\ast}} & X_{n}
}
\end{equation*}
gives the family $1_{X} = \{ 1_{X_{n}} \}$ the structure of a
simplicial presheaf on $\mathcal{C} / X$.  Note that the set
of all maps (\ref{eq 1}) can be identified with the collection of
ordinal number maps $\mathbf{n} \to \mathbf{k}$, and it follows that
there is an isomorphism of simplicial sets
\begin{equation*}
1_{X}(y) \cong \Delta^{k}
\end{equation*}
The canonical simplicial presheaf map $1_{X} \to \ast$ is therefore a
local weak equivalence. 

If $F$ is a presheaf on $\mathcal{C} / X$,
a presheaf map $f: 1_{X_{n}} \to F$ is completely determined by the
images under $f$ of the sections
\begin{equation*}
\xymatrix{
U \ar[r] \ar[d] & X_{n} \ar[d]^{1}  \\
X_{n} \ar[r]_{1} & X_{n}
}
\end{equation*}
so that
\begin{equation*}
\hom(1_{X_{n}},F) \cong \varprojlim_{x: U \to X_{n}}\ F_{n}(x)
\cong \Gamma_{\ast}F_{n},
\end{equation*}
where $F_{n}$ denotes the restriction of $F$ to the site $\mathcal{C}
/ X_{n}$. 

Thus, if $Z$ is a globally fibrant simplicial
presheaf on $\mathcal{C} / X$, there is a weak equivalence
\begin{equation*}
\Gamma_{\ast}Z = \mathbf{hom}(\ast,Z) \xrightarrow{\simeq}
\mathbf{hom}(1_{X},Z),
\end{equation*}
where the function space $\mathbf{hom}(1_{X},Z)$ is a homotopy inverse
limit of the simplicial sets $\Gamma_{\ast}Z_{n}$, computed on the
respective sites $\mathcal{C} / X_{n}$. This is, effectively,
an old observation --- see \cite{J1}.
\end{example}

\begin{example}\label{ex 4}
Suppose that $J$ is a small category, and identify $J$ with a constant
presheaf of categories on $\mathcal{C}$. The category $\mathcal{C}
/ J$ has, for objects, all pairs $(U,x)$ where $U$ is an
object of $\mathcal{C}$ and $x$ is an object of $J$, since $\Ob(J)$ is
a constant presheaf. The morphisms $(\alpha,f): (U,x) \to (V,y)$
are pairs consisting of a morphism $\alpha: U \to J$ of $\mathcal{C}$
and a morphism $f: x \to y$ of $J$. In other words, $\mathcal{C}
/ J = \mathcal{C} \times J$, and it's easy to see that this
identification gives $\mathcal{C} \times J$ the product topology, with
$J$ discrete. A sheaf on $\mathcal{C} \times J$ (and hence on
$\mathcal{C} / J$) can therefore be identified with a
$J^{op}$-diagram of sheaves on $\mathcal{C}$. 

Write $X_{j}$ for the simplicial presheaf $X(\ ,j)$ on $\mathcal{C}$.
A weak equivalence $X \to Y$ of simplicial presheaves on $\mathcal{C}
\times J$ is a map which induces a weak equivalence weak equivalences
$X_{j} \to Y_{j}$ of simplicial presheaves on $\mathcal{C}$ for all $j
\in J$. This can be proven directly by using the observation that $F$
is a sheaf on $\mathcal{C} \times J$ if and only if each $F_{j}$ is a
sheaf, or by using Lemma \ref{lem 14} below. It is also plain that a
map $A \to B$ of simplicial presheaves on $\mathcal{C} \times J$ is a
cofibration if and only if each map $A_{j} \to B_{j}$ is a cofibration
of simplicial presheaves on $\mathcal{C}$. The model structure for
simplicial presheaves on $\mathcal{C} \times J$ therefore coincides
with one of the standard model structures (due to Bousfield-Kan
\cite{BK}) for $J^{op}$-diagrams in the category $s\Pre(\mathcal{C})$
of simplicial presheaves on $\mathcal{C}$.

If $X$ is a globally
fibrant simplicial presheaf on $\mathcal{C} \times J^{op}$, then $X
\to \ast$ has the right lifting property with respect to all trivial
cofibrations $A \to B$ of $J^{op}$-diagrams of simplicial sets,
interpreted as trivial cofibrations of $(\mathcal{C} \times
J)$-presheaves which are constant in the $\mathcal{C}$
direction.  Then the global sections simplicial set $\Gamma_{\ast}X
= \varprojlim X$ can be written as an inverse limit
\begin{equation*}
\Gamma_{\ast}X = \varprojlim_{j} \Gamma_{\ast} X_{j}
\end{equation*}
of the global sections of the $\mathcal{C}$-presheaves $X_{j}$. The
indicated lifting property for $X$ means that the $J^{op}$-diagram $j
\mapsto \Gamma_{\ast}X_{j}$ of simplicial sets has the right lifting
property with respect to all trivial cofibrations of $J^{op}$
diagrams. It follows that $\Gamma_{\ast}X$ is the homotopy inverse
limit of the $J^{op}$-diagram $\Gamma_{\ast}X_{j}$.

Observe that the functor $Y \mapsto Y_{i}$ preserves weak
equivalences, and has a left adjoint defined by $A \mapsto A \times
\hom_{J^{op}}(\ ,i)$. This adjoint preserves trivial cofibrations, so
that $Y \mapsto Y_{i}$ preserves global fibrations. In particular, if
$X$ is a globally fibrant simplicial presheaf on $\mathcal{C} \times
J$, then all $X_{j}$ are globally fibrant simplicial presheaves on
$\mathcal{C}$.
\end{example}

\section{The fibred site for a presheaf}

Suppose that $X$ is a presheaf on $\mathcal{C}$. The corresponding
category $\mathcal{C} / X$ has as objects all pairs $(U,x)$
with $x \in X(U)$. The morphisms $(U,x) \to (V,y)$ of $\mathcal{C}
/ X$ consist of morphisms $\alpha: U \to V$ of $\mathcal{C}$
such that $\alpha^{\ast}(y) = x$. Such morphisms can be identified
with commutative diagrams of presheaf morphisms
\begin{equation*}
\xymatrix@C=5pt@R=10pt{
U \ar[rr]^{\alpha} \ar[dr]_{x} && V \ar[dl]^{y} \\
& X
}
\end{equation*}

Suppose that $\pi: Y \to X$ is a map of presheaves on $\mathcal{C}$, and
that $x: U \to X$ is an object of $\mathcal{C} / X$. Then $Y$
represents a presheaf $\pi_{\ast}$ on $\mathcal{C} / X$ by
setting $\pi_{\ast}(x)$ to be the set of sections
\begin{equation*}
\xymatrix{
& Y \ar[d]^{\pi} \\
U \ar[r]_{x} \ar[ur]^{\sigma} & X
}
\end{equation*}
of $\pi$ over $x$.

Conversely, if $F$ is a presheaf on $\mathcal{C} / X$, define
\begin{equation*}
F_{\ast}(U) = \bigsqcup_{x: U \to X}\ F(x).
\end{equation*}
Then any map $\alpha: U \to V$ in $\mathcal{C}$ defines a function
$\alpha^{\ast}: F_{\ast}(V) \to F_{\ast}(U)$, which is the unique
function making the diagrams
\begin{equation*}
\xymatrix{
F(y) \ar[r] \ar[d] & \bigsqcup_{y: V \to X}\ F(y) \ar[d]^{\alpha_{\ast}} \\
F(y\cdot \alpha) \ar[r] & \bigsqcup_{x: U \to X}\ F(x)
}
\end{equation*}
commute, where the horizontal functions are canonical. There is a
canonical function $\pi_{F}: F_{\ast}(U) \to X(U)$ which sends the
summand $F(x)$ to the section $x: U \to X$ in $X(U)$, and all such functions
form the components of a presheaf map $\pi_{F}: F_{\ast} \to X$.

The assignments $\pi \mapsto \pi_{\ast}$ and $F \mapsto \pi_{F}$ are
functorial, and define an equivalence of categories
\begin{equation*}
\Pre(\mathcal{C}) / X \simeq \Pre(\mathcal{C} / X).
\end{equation*}
Note that this equivalence specializes to an equivalence of presheaves
on $\mathcal{C} / U$ with morphisms of presheaves $Z \to U$
for each object $U$ of $\mathcal{C}$.

A presheaf morphism $\alpha: X \to Y$ induces a functor $\alpha: \mathcal{C}
/ X \to \mathcal{C} / Y$ by composition with
$\alpha$: the object $x: U \to X$ maps to the composite
\begin{equation*}
U \xrightarrow{x} X \xrightarrow{\alpha} Y.
\end{equation*}
Suppose that $F$ is a presheaf defined on $\mathcal{C} /
Y$. Then composition with $\alpha$ determines a presheaf $F\alpha$ on
$\mathcal{C} / X$, and it is easily seen that there is a
pullback diagram
\begin{equation*}
\xymatrix{
(F\alpha)_{\ast} \ar[r] \ar[d]_{\pi} & F_{\ast} \ar[d]^{\pi} \\
X \ar[r]_{\alpha} & Y
}
\end{equation*}

Any object $x: U \to X$ determines a functor $\phi_{x} = \phi_{U,x}:
\mathcal{C} / U \to \mathcal{C} / X$. As the
notation suggests, if $F$ is a presheaf on $\mathcal{C} / X$,
then the induced presheaf $F_{x} = F_{U,x}$ on $\mathcal{C} /
U$ is defined by composition with $x$. It follows that there is a
pullback diagram
\begin{equation}\label{eq 2} 
\xymatrix{
(F_{x})_{\ast} \ar[r] \ar[d]_{\pi} & F_{\ast} \ar[d]^{\pi} \\
U \ar[r]_{x} & X
}
\end{equation}
in the category of presheaves on $\mathcal{C}$.

Note that a presheaf map $\pi: Y \to X$ represents a sheaf on
$\mathcal{C} / X$ if and only if all presheaves
$(\pi_{\ast})_{U,x}$ of sections are sheaves on $\mathcal{C}
/ U$.
Equivalently, the map $\pi: Y \to X$ represents a sheaf on $\mathcal{C}
/ X$ if and only if, given a section $x \in X(U)$ and a
compatible family of sections
\begin{equation*}
\xymatrix{
& Y \ar[d] \\
U_{i} \ar[r]_{x_{i}} \ar[ur]^{\sigma_{i}} & X
}
\end{equation*}
defined over the restrictions $x_{i}$ of $x$ along some covering family $U_{i} \to U$, there is a unique section
\begin{equation*}
\xymatrix{
& Y \ar[d] \\
U \ar[r]_{x} \ar[ur]^{\sigma} & X
}
\end{equation*}
which restricts to all $\sigma_{i}$.

\begin{lemma}
The collection of all presheaf maps $\pi: Y \to X$ which represent
sheaves on $\mathcal{C} / X$ is stable under base change.
\end{lemma}

The statement of the Lemma (which is easy to prove) means that, given
a pullback square
\begin{equation*}
\xymatrix{
Z \times_{X} Y \ar[r] \ar[d]_{\pi_{\ast}} & Y \ar[d]^{\pi} \\
Z \ar[r] & X
}
\end{equation*}
if $\pi$ represents a sheaf on $\mathcal{C} / X$, then
$\pi_{\ast}$ represents a sheaf on $\mathcal{C} / Z$.

\begin{lemma}
A map $\pi: Y \to X$ represents a sheaf on $\mathcal{C} / X$
if and only if in all pullback diagrams
\begin{equation*}
\xymatrix{
U \times_{X} Y \ar[r] \ar[d]_{\pi_{\ast}} & Y \ar[d]^{\pi} \\
U \ar[r]_{x} & X
}
\end{equation*}
arising from sections $x \in X(U)$, $U \in \mathcal{C}$, the map
$\pi_{\ast}$ represents a sheaf on $\mathcal{C} / U$.
\end{lemma}

\begin{example}
Suppose that $G$ is a sheaf on the site $\mathcal{C}$, and let $q:
\mathcal{C} / U \to \mathcal{C}$ denote the canonical
forgetful functor. Then the composite $Gq$ is a sheaf on
$\mathcal{C} / U$. Explicitly, if $x: V \to U$ is an object
of $\mathcal{C} / U$, then $Gq(x) = G(V)$. It follows that 
\begin{equation*}
(Gq)_{\ast}(V) = \bigsqcup_{x: V \to U}\ G(V) = U(V) \times G(V),
\end{equation*}
and the canonical map $(Gq)_{\ast}(V) \to U(V)$ is just the projection
onto $U(V)$. In other words, $(Gq)_{\ast} \cong G \times U$, and the
canonical map $\pi$ is the projection $G \times U \to U$. This object
represents a sheaf on $\mathcal{C}/ U$ in the sense described
above, as one can check directly, but the product $G \times U$ need not
be a sheaf on the site $\mathcal{C}$.
\end{example}

We do, however, have the following:

\begin{lemma}
Suppose that $X$ is a sheaf and that $\pi: Y \to X$ is a presheaf map.
Then $\pi$ represents a sheaf on $\mathcal{C} / X$ if and
only if $Y$ is a sheaf.
\end{lemma}

\begin{proof}
The map $\pi: Y \to X$ represents a sheaf on $\mathcal{C}
/ X$ if and only if, given a section $x \in X(U)$ and a
compatible family of sections
\begin{equation*}
\xymatrix{
& Y \ar[d] \\
U_{i} \ar[r]_{x_{i}} \ar[ur]^{\sigma_{i}} & X
}
\end{equation*}
defined over the restrictions $x_{i}$ of $x$ along some covering family $U_{i} \to U$, there is a unique section
\begin{equation*}
\xymatrix{
& Y \ar[d] \\
U \ar[r]_{x} \ar[ur]^{\sigma} & X
}
\end{equation*}
which restricts to all $\sigma_{i}$.

If $Y$ is a sheaf on $\mathcal{C}$, then there is a unique element
$\sigma: U \to Y$ which restricts to all $\sigma_{i}$. Since $X$ is a
sheaf, $\pi(\sigma) = x$.

Suppose that $\pi$ represents a sheaf on $\mathcal{C} / X$,
and let $\sigma_{i}: U_{i} \to Y$ be a compatible family of elements
defined on a covering $U_{i} \to U$ of $U$. Then $\pi(\sigma_{i}) =
x_{i}$ and the $x_{i}$ uniquely determine a section $x: U \to X$ since
$X$ is a sheaf on $\mathcal{C}$. But then a lifting $\sigma: U \to Y$
of $x$ exists and is uniquely determined since $\pi$ represents a
sheaf on $\mathcal{C} / X$. Any such lifting $\sigma$
extending the $\sigma_{i}$ must map to $x$, since $X$ is a sheaf.
\end{proof}

Say that a map 
\begin{equation*}
\xymatrix@C=5pt@R=10pt{
Z \ar[rr]^{f} \ar[dr] && W \ar[dl] \\
& X
}
\end{equation*}
of simplicial presheaves over $X$ is a local weak equivalence fibred
over $X$ if it represents a local weak equivalence of simplicial
presheaves on $\mathcal{C} / X$.

\begin{lemma}\label{lem 9}
Suppose that $X$ is a presheaf on $\mathcal{C}$.
Suppose that 
\begin{equation*}
\xymatrix@C=5pt@R=10pt{
Z \ar[rr]^{f} \ar[dr] && W \ar[dl] \\
& X
}
\end{equation*}
is a commutative diagram of simplicial presheaves. Then $f$ represents
a local weak equivalence of simplicial presheaves on $\mathcal{C}
/ X$ if and only if the map $Z \to W$ is a local weak
equivalence of simplicial presheaves on $\mathcal{C}$.
\end{lemma}

\begin{proof}
Recall that if $Z \to X$ is a map of simplicial presheaves, then the
presheaf that it represents on $\mathcal{C} / X$ associates
to $x: V \to X$, the fibre $Z_{x}$ over $x$ for the simplicial set map
$Z(V) \to X(V)$. Certainly, $Z(V) = \sqcup_{x \in X(V)} X_{x}$ and
it's clear that $\Ex^{\infty}Z(V) = \sqcup_{x \in X(V)}
\Ex^{\infty}Z_{x}$. In particular, the map $Z \to \Ex^{\infty}Z$
fibres over $X$, and represents a sectionwise weak equivalence of
simplicial presheaves on $\mathcal{C}/ X$. The map $Z \to
\Ex^{\infty}Z$ is also a sectionwise weak equivalence
of simplicial presheaves on $\mathcal{C}$. It suffices, therefore, to
assume that $Z$ and $W$ are presheaves of Kan complexes on $\mathcal{C}$.

In that case, the map $f$ has the standard factorization
\begin{equation*}
\xymatrix{
Z \ar[r]^{j} \ar[dr]_{f_{\ast}} & T \ar[d]^{p} \\
& W
}
\end{equation*}
where $p$ is a sectionwise Kan fibration and $j$ is left inverse to a
sectionwise trivial Kan fibration. Furthermore, this factorization is
fibred over $X$ and has the same properties in each fibre. In
particular $j$ represents a sectionwise trivial map of simplicial
presheaves on $\mathcal{C} / X$. It suffices, therefore, to
assume that $f$ is a sectionwise Kan fibration between presheaves of
Kan complexes, and show that it is locally trivial for the site
$\mathcal{C}$ if and only if it is locally trivial for the site
$\mathcal{C} / X$.

Suppose given a commutative diagram
\begin{equation*}
\xymatrix{
\partial\Delta^{n} \ar[r] \ar[d]  & Z_{x} \ar[d]^{f} \\
\Delta^{n} \ar[r] & W_{x}
}
\end{equation*}
where $x: V \to X$ is an object of $\mathcal{C} / X$. Then
there is a covering family $\phi_{i}: V_{i} \to V$ for which the
displayed liftings exist in the diagram
\begin{equation*}
\xymatrix{ 
\partial\Delta^{n} \ar[r] \ar[d]  & Z_{x} \ar[r] 
& Z(V) \ar[r] & Z(V_{i}) \ar[d]^{f} \\
\Delta^{n} \ar[r] \ar[urrr]^{\sigma_{i}} 
& W_{x} \ar[r] & W(V) \ar[r]
& W(V_{i})
}
\end{equation*}
But then $\sigma_{i}$ factors through the summand
$Z_{\phi_{i}^{\ast}(x)}$, since its image in $W(V_{i})$
factors through the summand $W_{\phi_{i}^{\ast}(x)}$.

Conversely, suppose given a diagram
\begin{equation*}
\xymatrix{
\partial\Delta^{n} \ar[r]^{\alpha} \ar[d] & Z(V) \ar[d]^{f} \\
\Delta^{n} \ar[r]_-{\beta} & W(V)
}
\end{equation*}
Then $\beta$ factors through a summand $W_{x}$ for some $x: V \to X$
in $\mathcal{C}$ since $\Delta^{n}$ is connected, and it follows that
$\alpha$ factors through the summand $Z_{x}$. Thus if the lifting
problem can be solved locally over $\mathcal{C} / X$ it can
be solved locally over $\mathcal{C}$. It follows that if $f$
represents a local trivial fibration on the site $\mathcal{C}
/ X$, then $f$ is a local trivial fibration on $\mathcal{C}$.
\end{proof}

\begin{corollary}
Suppose that $X$ is a presheaf on $\mathcal{C}$. The model
structure on the category $s\Pre(\mathcal{C}) / X$ which arises
from the topology on $\mathcal{C} / X$ is induced
from the model structure on the category
$s\Pre(\mathcal{C})$ of simplicial presheaves. In particular, a map
\begin{equation*}
\xymatrix@C=5pt@R=10pt{
Z \ar[rr]^{f} \ar[dr] && W \ar[dl] \\
& X
}
\end{equation*}
is a weak equivalence (respectively cofibration, fibration) if and
only if the map $f: Z \to W$ is a weak equivalence (respectively
cofibration, global fibration) of simplicial presheaves on
$\mathcal{C}$.
\end{corollary}

\begin{proof}
The statement about cofibrations amounts to the observation that
cofibrations are defined fibrewise. The statement for weak
equivalences is Lemma \ref{lem 9}, and then the fibration statement
is a formal consequence.
\end{proof}

\begin{lemma}\label{lem 11}
\begin{itemize}
\item[1)]
Suppose that $\alpha: X' \to X$ is a morphism of presheaves. Then the
functor $s\Pre(\mathcal{C}) / X \to s\Pre(\mathcal{C})
/ X'$ defined by pullback preserves weak equivalences.
\item[2)]
Suppose that $X$ is a presheaf on $\mathcal{C}$.
Then a map
\begin{equation*}
\xymatrix@C=5pt@R=10pt{
Z \ar[rr]^{f} \ar[dr] && W \ar[dl] \\
& X
}
\end{equation*}
of simplicial presheaves over $X$ represents a local weak equivalence
on $\mathcal{C} / X$ if and only if all pullbacks
\begin{equation*}
\xymatrix@C=5pt{
U \times_{X} Z \ar[rr]^{f} \ar[dr] && U \times_{X} W \ar[dl] \\
& U
}
\end{equation*}
over sections $x: U \to X$, $U \in \mathcal{C}$ represent local weak
equivalences on $\mathcal{C} / U$.
\end{itemize}
\end{lemma}

\begin{proof}
Statement 2) implies statement 1). We shall prove statement 2).

Recall that 
\begin{equation*}
Z(V) = \bigsqcup_{x \in X(V)}\ Z_{x}
\end{equation*}
 and that 
\begin{equation*}
U \times_{X} Z(V) = \bigsqcup_{\phi: V \to U}\ Z_{x\phi}. 
\end{equation*}
Once again,
the $\Ex^{\infty}$ construction is performed fibrewise, so it suffices
to assume that $Z$ and $W$ are presheaves of Kan complexes. The
canonical replacement of a map by a fibration is also a fibrewise
construction, so it suffices to assume that $f$ is a Kan fibration in
each section, and hence in each fibre. But then
$f$ has the local right lifting property with respect
to all inclusions $\partial\Delta^{n} \subset \Delta^{n}$ if and only
if all pullbacks of $f$ along sections $x: U \to X$ have the same
local right lifting property, by the argument that appears in the
proof of Lemma \ref{lem 9}.    
\end{proof}

Note that statement 1) of Lemma \ref{lem 11} is not true if $X'$ and
$X$ are replaced by simplicial presheaves. One can see counterexamples
easily in ordinary simplicial sets.

\section{Constructions for presheaves of categories}

Suppose that $A$ is a presheaf of categories. An {\it enriched
diagram} $X$ on $A$ consists of set-valued functors $X(U): A(U) \to
\mathbf{Set}$ defined by $x \mapsto A(U)_{x}$, one for each $U \in
\mathcal{C}$, such that each morphism $\phi: V \to U$ of $\mathcal{C}$
induces functions $\phi^{\ast}: X(U)_{x} \to X(V)_{\phi^{\ast}(x)}$
and all diagrams
\begin{equation*}
\xymatrix{
X(U)_{x} \ar[r]^{\alpha_{\ast}} \ar[d]_{\phi^{\ast}} 
& X(U)_{y} \ar[d]^{\phi^{\ast}} \\
X(V)_{\phi^{\ast}(x)} \ar[r]_{(\phi^{\ast}(\alpha))_{\ast}} 
& X(V)_{\phi^{\ast}(y)}
}
\end{equation*} 
commute, where $\alpha: x \to y$ is a morphism of $A(U)$.

Let $F$ be a presheaf on the fibred site $\mathcal{C} /
A$. Then $F$ assigns a set $F(U)_{x} = F(U,x)$ to each object $x: U
\to \Ob(A)$. Every morphism $\gamma: x \to y$ in $A(U)$ determines a
morphism $(1,\gamma): (U,x) \to (U,y)$, and hence induces a function
$(1,\gamma)^{\ast}: F(U)_{y} \to F(U)_{x}$. In particular, $F$
determines a functor $F(U): A(U)^{op} \to \mathbf{Set}$. Any
morphism $\alpha: V \to U$ of $\mathcal{C}$ induces a morphism
$(\alpha,1): (V,\alpha^{\ast}(x)) \to (U,x)$ in $\mathcal{C}
/ A$, and hence induces a function $\alpha^{\ast}: F(U)_{x}
\to F(V)_{\alpha^{\ast}(x)}$. If $\alpha: V \to U$ is a morphism of
$\mathcal{C}$ and $\gamma: x \to y$ is a morphism of $A(U)$ then the
diagram
\begin{equation*}
\xymatrix{
(V,\alpha^{\ast}(x)) 
\ar[r]^{(1,\alpha^{\ast}(\gamma))} \ar[d]_{(\alpha,1)} 
& (V,\alpha^{\ast}(y)) \ar[d]^{(\alpha,1)} \\
(U,x) \ar[r]_{(1,\gamma)} & (U,y)
}
\end{equation*} 
commutes in $\mathcal{C}/ A$, so that the diagram
\begin{equation*}
\xymatrix{
F(U)_{y} \ar[r]^{(1,\gamma)^{\ast}} \ar[d]_{(\alpha,1)^{\ast}} 
& F(U)_{x} \ar[d]^{(\alpha,1)^{\ast}} \\
F(V)_{\alpha^{\ast}(y)} \ar[r]_{(1,\alpha^{\ast}(\gamma))^{\ast}} 
& F(U)_{\alpha^{\ast}(x)}
}
\end{equation*}
commutes. In other words, $F$ defines an enriched diagram $F$ on the
presheaf of categories $A^{op}$.

Suppose that $G$ is an enriched diagram on the presheaf of categories
$A^{op}$. Write $G(U,x) = G(U)_{x}$ for each object $(U,x)$ of
$\mathcal{C} / A$. Let $(\alpha,\gamma): (V,y) \to (U,x)$ be
a morphism of $\mathcal{C} / A$. Then $(\alpha,\gamma)$ has
a factorization
\begin{equation*}
\xymatrix{
(V,y) \ar[r]^{(1,\gamma)} \ar[dr]_{(\alpha,\gamma)} 
& (V,\alpha^{\ast}(x)) \ar[d]^{(\alpha,1)} \\
& (U,x)
}
\end{equation*}
Associate to $(\alpha,\gamma)$ the composite
\begin{equation*}
G(U)_{x} \xrightarrow{\alpha^{\ast}} G(V)_{\alpha^{\ast}(x)} 
\xrightarrow{\gamma^{\ast}} G(V)_{y}.
\end{equation*}
If $(\beta,\omega): (W,z) \to (V,y)$ is another choice of
morphism of $\mathcal{C} / A$, there is a commutative
diagram
\begin{equation*}
\xymatrix{
(W,z) \ar[r]^{(1,\omega)} \ar[dr]_{(\beta,\omega)} 
& (W,\beta^{\ast}(y)) \ar[d]^{(\beta,1)} \ar[r]^{(1,\beta^{\ast}(\gamma))} 
& (W,\beta^{\ast}\alpha^{\ast}(x)) \ar[d]^{(\beta,1)} \\
& (V,y) \ar[r]^{(1,\gamma)} \ar[dr]_{(\alpha,\gamma)} 
& (V,\alpha^{\ast}(x)) \ar[d]^{(\alpha,1)} \\
&& (U,x)
}
\end{equation*}
It follows that the assignment $(U,x) \mapsto G(U)_{x}$ defines a
presheaf on the category $\mathcal{C} / A$.

We have proved the following:

\begin{lemma}\label{lem 12}
Suppose that $A$ is a presheaf of categories on a small site
$\mathcal{C}$. Then the category $\Pre(\mathcal{C} / A)$ is
equivalent to the category of enriched diagrams on the presheaf of
categories $A^{op}$.
\end{lemma}

Note that a presheaf $F$ on $\mathcal{C} / A$ consists of a
presheaf of objects $F_{0} \to \Ob(A)$ on $\mathcal{C} /
\Ob(A)$, with extra structure.  

There is a canonical functor $\psi:
\Ob(A) \to A$, and the assignment $F \mapsto F_{0}$ coincides with the
restriction functor
\begin{equation*}
\psi_{\ast}: \Pre(\mathcal{C} / A) \to \Pre(\mathcal{C}
/ \Ob(A))
\end{equation*}
which is defined by composition with the canonical functor $\psi$,
under the equivalence
\begin{equation*}
\Pre(\mathcal{C} / \Ob(A) \simeq \Pre(\mathcal{C}) / \Ob(A)
\end{equation*}
of the last section.

The object $(U,x)$ of the site
$\mathcal{C} / A$ determines a functor
\begin{equation*}
\phi_{U,x}: \mathcal{C} / U \to \mathcal{C} / A
\end{equation*}
which sends an object $\phi: V \to U$ to the object
$(W,\phi^{\ast}x)$. This functor sends the morphism
\begin{equation*}
\xymatrix@C=5pt@R=10pt{
V \ar[dr]_{\phi} \ar[rr]^{\alpha} && V' \ar[dl]^{\phi'} \\
& U 
}
\end{equation*}
to the morphism
\begin{equation*}
(W,\phi^{\ast}(x)) \xrightarrow{(\alpha,1)} (W',(\phi')^{\ast}(x)).
\end{equation*}
When $F$ is a presheaf on $\mathcal{C} / A$,
write $F_{U,x}$ for the presheaf on $\mathcal{C} / U$ which
is defined by composition with $\phi_{U,x}$ in the sense that 
\begin{equation*}
F_{U,x} = F \cdot \phi_{U,x}.
\end{equation*}

Now here are some observations:

\begin{lemma}
\begin{itemize}
\item[1)]
A presheaf $F$ on $\mathcal{C}/ A$ is a sheaf if and only if
all restricted presheaves $F_{U,x}$ are sheaves on $\mathcal{C}
/ U$.
\item[2)]
A presheaf $F$ on $\mathcal{C} / A$ is a sheaf if and only if
it restricts to a sheaf $F_{0} \to \Ob(A)$ on $\mathcal{C} / \Ob(A)$.
\item[3)]
The restrictions $F \mapsto F_{U,x}$ and $F \mapsto F_{0}$ commute with
the associated sheaf functor on $\mathcal{C} / U$ and
$\mathcal{C} / \Ob(A)$ respectively, up to natural isomorphism.
\end{itemize}
\end{lemma}

In the same way, a simplicial presheaf $X$ on the fibred site
$\mathcal{C} / A$ consists of a simplicial presheaf of
objects $X_{0} \to \Ob(A)$ over the presheaf $\Ob(A)$ with extra
structure. The restriction functor
\begin{equation*}
\psi_{\ast}: s\Pre(\mathcal{C} / A) \to s\Pre(\mathcal{C} /
\Ob(A))
\end{equation*}
can be identified up to equivalence with the object functor
\begin{equation*}
s\Pre(\mathcal{C} / A) \to s\Pre(\mathcal{C}) /
\Ob(A)
\end{equation*}
which takes a simplicial presheaf $X$ (aka. enriched diagram in
simplicial sets) to the simplicial presheaf of objects $X_{0} \to
\Ob(A)$ over $\Ob(A)$.

\begin{lemma}\label{lem 14}
The object functor $\psi_{\ast}: s\Pre(\mathcal{C} / A) \to
s\Pre(\mathcal{C}) / \Ob(A)$ preserves and reflects local
weak equivalences.
\end{lemma}

\begin{proof}
We show that a map $f: X \to Y$ of simplicial presheaves on
$\mathcal{C} / A$ is a local weak equivalence if and only if
the induced map $X_{0} \to Y_{0}$ is a local weak equivalence of
simplicial presheaves. The statement of the result is a generalization
of Lemma \ref{lem 9}, and the proof involves the same ideas.
This works because the topology on $\mathcal{C} / A$ only
involves the topology on $\mathcal{C} / \Ob(A)$. 

The forgetful functor preserves sectionwise weak equivalences.
A map $f: X \to Y$ is a local weak equivalence if and only if the
induced map $\Ex^{\infty}X \to \Ex^{\infty}Y$ is a local weak
equivalence. The canonical map $j: X \to \Ex^{\infty}X$ is a sectionwise
weak equivalence. The $\Ex^{\infty}$ construction and the associated
sectionwise equivalence are preserved by the forgetful functor. Thus
it suffices to assume that $X$ and $Y$ are presheaves of Kan complexes.

In that case the map $f: X \to Y$ has a standard factorization
\begin{equation*}
\xymatrix{
X \ar[r]^{i} \ar[dr]_{f} & Z \ar[d]^{p} \\
& Y
}
\end{equation*}
where $p$ is a Kan fibration in each section and $i$ is right inverse
to a sectionwise trivial Kan fibration. It therefore suffices to
assume that $f$ is a Kan fibration in each section, and show that $f$
is a local trivial fibration if and only if the induced map $f_{0}:
X_{0} \to Y_{0}$ is a local trivial fibration. But this is now clear:
the argument is finished as in the proof of Lemma \ref{lem 9}.
\end{proof}

The object functor $s\Pre(\mathcal{C} / A) \to
s\Pre(\mathcal{C}) / \Ob(A)$ also preserves and reflects
monomorphisms. 

The restriction functor $\psi_{\ast}$ has a left adjoint 
\begin{equation*}
\psi^{\ast}: s\Pre(\mathcal{C}
/ \Ob(A)) \to s\Pre(\mathcal{C} / A).
\end{equation*}
which is defined by left Kan extension along the inclusion $\psi:
\Ob(A) \to A$, and we identify this with a left adjoint
\begin{equation*}
\psi^{\ast}: s\Pre(\mathcal{C})
/ \Ob(A) \to s\Pre(\mathcal{C} / A).
\end{equation*}
for the object functor. For a fixed simplicial presheaf $X \to
 \Ob(A)$ over $\Ob(A)$, the map $\psi^{\ast}X_{0} \to \Ob(A)$ can be
 identified with the composite
\begin{equation*}
X \times_{\Ob(A)} \Mor(A) \to \Mor(A) \xrightarrow{t} \Ob(A) 
\end{equation*}
where $t$ is the target map, and both the indicated pullback and the
projection are defined by the pullback diagram
\begin{equation*}
\xymatrix{
X \times_{\Ob(A)} \Mor(A) \ar[r] \ar[d] & \Mor(A) \ar[d]^{s} \\
X \ar[r] & \Ob(A)
}
\end{equation*}
Here, $s$ is the source map. The map $s$ is a local fibration since
$\Mor(A)$ and $\Ob(A)$ are simplicial presheaves which are constant in
the simplicial direction. It follows that the indicated
pullback is a homotopy cartesian diagram of simplicial presheaves. It
follows that the functor which sends the simplicial presheaf map $X
\to \Ob(A)$ to $(\psi^{\ast}X)_{0} = X \times_{\Ob(A)} \Mor(A)$
preserves local weak equivalences. It also preserves
cofibrations. This suffices for a proof of the following:

\begin{lemma}\label{lem 15}
The object functor $\psi_{\ast}: s\Pre(\mathcal{C} / A) \to
s\Pre(\mathcal{C}) / \Ob(A)$ defined by sending $X$ to the
map $X_{0} \to \Ob(A)$ preserves global fibrations.
\end{lemma}

In particular, a global fibration $X \to Y$ in $s\Pre(\mathcal{C}
/ A)$ consists of a global $X_{0} \to Y_{0}$ over $\Ob(A)$
which is $A$-equivariant in an enriched sense.

For a fixed simplicial presheaf (or enriched functor) $X$ on
$\mathcal{C} / A$, applying the homotopy colimit functor in
each section gives a simplicial presheaf $\hocolim_{A^{op}}\ X$ and a
canonical simplicial presheaf map $\pi: \hocolim_{A}\ X \to BA^{op}$.
This assignment is plainly functorial in $X$.

\begin{lemma}\label{lem 16}
The homotopy colimit functor $s\Pre(\mathcal{C} / A) \to
s\Pre(\mathcal{C}) / BA^{op}$ preserves weak equivalences.
\end{lemma}

\begin{proof}
Note first of all that $\Ob(A^{op}) = \Ob(A)$.

There is a presheaf $Mor_{n}(A^{op})$ which consists of strings of
arrows of length $n$ in the presheaf of categories $A^{op}$, and
$\hocolim_{A^{op}}\ X$ is the diagonal of a bisimplicial sheaf which
is given by the object $X \times_{\Ob(A)} \Mor_{n}A^{op}$ in horizontal
degree $n$. Here, the map $s_{0}: \Mor_{n}(A^{op}) \to \Ob(A)$ is defined
by picking out the first object in the string, and is a local
fibration. It follows that the pullback diagram of simplicial presheaf
maps
\begin{equation*}
\xymatrix{
X \times_{\Ob(A)} \Mor_{n}(A^{op}) \ar[r] \ar[d] 
& \Mor_{n}(A^{op}) \ar[d]^{s_{0}}  \\
X \ar[r] & \Ob(A)
}
\end{equation*}
is homotopy cartesian, so that any local weak equivalence $X \to Y$ over
$\Ob(A)$ induces a local weak equivalence 
\begin{equation*}
X \times_{\Ob(A)} \Mor_{n}(A^{op}) \to Y \times_{\Ob(A)} \Mor_{n}(A^{op}).
\end{equation*}
This is true in  all horizontal degrees $n$, and so the map
\begin{equation*}
\hocolim_{A^{op}}\ X \to \hocolim_{A^{op}}\ Y
\end{equation*}
is a local weak equivalence.
\end{proof}

The model structure that we have been using so far on 
$s\Pre(\mathcal{C} / A)$ is the natural ``injective''
structure, for which a map $f: X \to Y$ of enriched diagrams in
simplicial presheaves is a weak equivalence (respectively cofibration)
if and only if the induced map
\begin{equation*}
\xymatrix@C=5pt@R=10pt{
X_{0} \ar[rr]^{f_{0}} \ar[dr] && Y_{0} \ar[dl] \\
& \Ob(A)
}
\end{equation*}
is a weak equivalence (respectively cofibration) of
$s\Pre(\mathcal{C}) / \Ob(A)$. There is also a projective
structure on $s\Pre(\mathcal{C} / A)$ which has the same weak
equivalences, but for which a map $f$ is a fibration if and only if
the induced diagram as above is a fibration of $s\Pre(\mathcal{C})
/ \Ob(A)$. Say that such a map is a projective fibration, and
say that a projective cofibration is a map which has the left lifting
property with respect to all maps $p: X \to Y$ which are
simultaneously projective fibrations and local weak equivalences.

\begin{lemma}\label{lem 17}
The category $s\Pre(\mathcal{C} / A)$ of enriched diagrams on
$A^{op}$, together with the local weak equivalences, projective
fibrations and projective cofibrations as defined above, satisfies the
axioms for a closed model category.
\end{lemma} 

\begin{proof}
A map $p: X \to Y$ is a projective fibration (respectively trivial
projective fibration) if and only if it has the right lifting property
with respect to all maps $i_{\ast}: \psi^{\ast}A \to \psi^{\ast}B$
where $i:A \to B$ is a trivial cofibration (respectively cofibration)
over $\Ob(A)$. We have already seen that the functor $A \mapsto
\psi^{\ast}A$ preserves local weak equivalences. The factorization
axiom is now an easy consequence of these observations, along with the
standard fact that the ``injective'' model structure for the category
of simplicial presheaves is cofibrantly generated. The lifting axiom
{\bf CM4} follows by a standard argument.
\end{proof}

Suppose that $\phi: A \to B$ is a functor of presheaves
of categories. Then precomposition with $\phi$ defines a restriction functor
\begin{equation*}
\phi_{\ast}: s\Pre(\mathcal{C} / B) \to s\Pre(\mathcal{C}
/ A).
\end{equation*}
In effect, an enriched diagram $X$ on $B$ taking values in simplicial
sets consists of contravariant simplicial set-valued functors $X(U):
B(U) \to \mathbf{S}$, $U \in \mathcal{C}$ which fit together along
morphisms of $\mathcal{C}$, and then $\phi_{\ast}X$ consists of the
composite functors
\begin{equation*}
A(U) \xrightarrow{\phi} B(U) \xrightarrow{X} \mathbf{S}.
\end{equation*}

The following result is a corollary of Lemma \ref{lem 11} and Lemma
\ref{lem 14}:

\begin{corollary}\label{cor 18}
The restriction functor $\phi_{\ast}$ preserves local weak
equivalences for any functor $\phi: A \to B$ of presheaves of
categories.
\end{corollary}

\begin{proof}
The object functor $\psi_{\ast}: s\Pre(\mathcal{C} / A) \to
s\Pre(\mathcal{C}) / \Ob(A)$ detects weak equivalences, and
there is a relation $\phi_{\ast}\psi_{\ast} =
\psi_{\ast}\phi_{\ast}$. The functor $\phi_{\ast}$ induced by the
object-level morphism $\phi: \Ob(A) \to \Ob(B)$ preserves weak
equivalences by Lemma \ref{lem 11}.
\end{proof}

Note that there is a pullback diagram of simplicial presheaves
\begin{equation*}
\xymatrix{
(\phi_{\ast}X)_{0} \ar[r] \ar[d] & X_{0} \ar[d] \\
\Ob(A) \ar[r]_{\phi} & \Ob(B)
}
\end{equation*}
The functor $X \mapsto \phi_{\ast}X$ preserves projective fibrations
almost by definition, and it follows from Corollary \ref{cor 18} that
$\phi_{\ast}$ preserves trivial projective fibrations. The functor
$\phi_{\ast}$ therefore determines a derived functor
\begin{equation*}
R\phi_{\ast}: \Ho(s\Pre(\mathcal{C} / B)) \to \Ho(s\Pre(\mathcal{C}
/ A))
\end{equation*}
which is defined by $R\phi_{\ast}(X) = \phi_{\ast}FX$, where the
trivial cofibration $j: X \to FX$ is a projective fibrant replacement for
$X$. 

The left adjoint
\begin{equation*}
\phi^{\ast}: s\Pre(\mathcal{C} / A) \to s\Pre(\mathcal{C}
/ B),
\end{equation*}
preserves projective
cofibrations and weak equivalences between projective cofibrant
objects, and therefore has an associated derived functor
\begin{equation*}
L\phi^{\ast}: \Ho(s\Pre(\mathcal{C} / A)) \to \Ho(s\Pre(\mathcal{C}
/ B)),
\end{equation*}
which is left adjoint to the derived functor $R\phi_{\ast}$.  The
derived functor $L\phi^{\ast}$ is defined by $L\phi^{\ast}(Y) =
\phi^{\ast}CY$, where the trivial projective fibration $p: CY \to Y$
is a projective cofibrant replacement for $Y$.

\section{Simplicial set constructions}

Suppose that $G$ is a groupoid and that $A: G \to
s\mathbf{Set}$ is a $G$-diagram in the category of simplicial
sets --- write $s\mathbf{Set}^{G}$ for the category of all such
objects. The diagram $A$ determines a canonical simplicial set map
$\hocolim_{G}\ A \to BG$, where $\hocolim_{G}\ A$ is
identified with the diagonal of the usual bisimplicial set.

In general, if $f: X \to BG$ is a simplicial set map, then $f$
can be identified with a set-valued functor $\sigma \mapsto
X_{\sigma}$ defined on the simplex category $\mathbf{\Delta}
/ BG$ of $BG$, where $X_{\sigma} = f^{-1}(\sigma)$
is the fibre over $\sigma$ for the function $X_{n} \to BG_{n}$ if
$\sigma$ is an $n$-simplex of $BG$. Note that a morphism
\begin{equation*}
\xymatrix@R=5pt{
\Delta^{m} \ar[dr]^{\theta^{\ast}\sigma} \ar[dd]_{\theta} & \\
& BG \\
\Delta^{n} \ar[ur]_{\sigma}
}
\end{equation*}
induces a function $X_{\sigma} \to X_{\theta^{\ast}(\sigma)}$ in the
obvious way. Suppose that $\sigma$ is the string 
\begin{equation*}
a_{0} \to a_{1} \to \dots \to a_{n}
\end{equation*}
of morphisms of $G$. Then a map 
\begin{equation*}
\xymatrix@C=5pt{
X \ar[rr]^-{g} \ar[dr]_{f} && \hocolim_{G}\ A \ar[dl]^{\pi} \\
& BG
}
\end{equation*}
of simplicial sets over $BG$ can be identified with a natural
transformation $g: X_{\sigma} \to (A_{a_{0}})_{n}$ over the simplex category
of $BG$; the naturality means that all diagrams
\begin{equation}\label{eq 3}
\xymatrix{
X_{\sigma} \ar[r]^{g} \ar[dd]_{\theta^{\ast}} & (A_{a_{0}})_{n} \ar[d]^{\theta^{\ast}} \\
& (A_{a_{0}})_{m} \ar[d]^{\theta_{\ast}} \\
X_{\theta^{\ast}(\sigma)} \ar[r]_{g} & (A_{a_{\theta(0)}})_{m}
}
\end{equation}
commute, where $\theta_{\ast}$ is induced by the map $a_{0} \to a_{\theta(0)}$
of $G$.

Suppose that $y$ is an object of $G$, let $f: X \to BG$ be a
simplicial set map, and define $\pb(X)_{y}$ by the pullback diagram
\begin{equation}\label{eq 4}
\xymatrix{
\pb(X)_{y} \ar[r] \ar[d] & X \ar[d]^{f} \\
B(G/ y) \ar[r] & BG
}
\end{equation}
where $B(G / y) \to BG$ is induced by the forgetful
functor $G / y \to G$. An $n$-simplex of
$\pb(X)_{y}$ consists of a triple
\begin{equation*}
(x,\sigma: a_{0} \to \dots \to a_{n},\alpha: a_{n} \to y), 
\end{equation*}
where $x \in X_{n}$, $f(x) = \sigma$ and $\alpha$ is a morphism of
$G$. Since $G$ is a groupoid, all morphisms in the string
$\sigma$ are invertible, and we can instead identify the $n$-simplex
of $\pb(X)_{y}$ displayed by the triple above, with a triple of the form
\begin{equation*}
(x,\sigma: a_{0} \to \dots \to a_{n},\gamma: a_{0} \to y). 
\end{equation*}
Of course, the assignment $y \to \pb(X)_{y}$ defines a functor $\pb(X): G \to
s\mathbf{Set}$. 

Observe that there is an inclusion
\begin{equation*}
c_{\sigma}: X_{\sigma} \to (\pb(X)_{a_{0}})_{n}
\end{equation*}
which is defined by sending $x$ to the triple $(x,\sigma,1: a_{0} \to
a_{0})$. It is not hard to show that diagrams of the form (\ref{eq 3})
commute for the list of functions $\{ c_{\sigma} \}$, and so these
functions define a natural map $\eta: X \to \hocolim_{G}\ \pb(X)$ of
simplicial sets over $BG$. 

There is a simplicial map $\epsilon_{y}: \pb(\hocolim_{G}\ A)_{y}
\to A_{y}$ which is defined on $n$-simplices by sending the triple
\begin{equation*}
(x \in A_{a_{0}},\sigma,\gamma: a_{0} \to y)
\end{equation*}
to the element $\gamma_{\ast}(x) \in (A_{y})_{n}$. This map is natural
in $y$ and in $A$, and therefore defines a natural map of
$G$-diagrams $\epsilon: \pb(\hocolim_{G}\ A) \to A$. It
is not difficult to show that the natural maps $\eta$ and $\epsilon$
satisfy the triangle identities, so that we have proved

\begin{lemma}
Suppose that $G$ is a groupoid.  Then the functor $\pb$ is left
adjoint to the homotopy colimit functor $\hocolim:
s\mathbf{Set}^{G} \to s\mathbf{Set}/ BG$.
\end{lemma}

\begin{lemma}\label{lem 20}
The canonical map $c: \hocolim_{G}\ \pb(X) \to X$ is a weak
equivalence, for all objects $f:X \to BG$ of the category
$s\mathbf{Set}/ BG$ of simplicial sets over $BG$.
\end{lemma}

\begin{proof}
The map $c$ is induced by a map of bisimplicial sets which is
specified in horizontal degree $n$ by the simplicial set map 
\begin{equation*}
\bigsqcup_{y_{0} \to \dots \to y_{n}}\ \pb(X)_{y_{0}} \to X.
\end{equation*}
which will also be denoted by $c$.
Note that $BG \cong \varinjlim_{y \in G} B(G/ y)$,
so that $X \cong \varinjlim_{y \in G} \pb(X)_{y}$. If $x \in X_{\sigma} =
f^{-1}(\sigma)$, where $\sigma$ is the $k$-simplex $z_{0} \to \dots
\to z_{k}$ of $BG$, then the preimage of $x$ under $c$
can be identified with a copy of 
$B(z_{k} / G)$, which is contractible. It follows
that the bisimplicial set map $c$ is a weak equivalence in each vertical
degree, and therefore induces a weak equivalence of associated
diagonals.
\end{proof}

\begin{corollary}\label{cor 21}
The map $\eta: X \to \hocolim_{G}\ \pb(X)$ is a weak equivalence.
\end{corollary}

\begin{proof}
The map $\hocolim_{G} \pb(X) \to X$ of Lemma \ref{lem 20} is a
left inverse for $\eta$.
\end{proof}

\begin{corollary}\label{cor 22}
The counit map $\epsilon: \pb(\hocolim_{G}\ A) \to A$ is a
weak equivalence for all $G$-diagrams $A$.
\end{corollary}

\begin{proof}
The induced map $\hocolim_{G}\pb(\hocolim_{G}\ A) \to
\hocolim_{G}\ A$ is a weak equivalence, by Corollary \ref{cor
21} together with the fact that $\eta$ and $\epsilon$ satisfy the
triangle identities. At the same time, all diagrams
\begin{equation*}
\xymatrix{
B_{y} \ar[r] \ar[d] & \hocolim_{G}\ B  \ar[d]^{\pi} \\
\Delta^{0} \ar[r]_{y} & BG
}
\end{equation*}
are homotopy cartesian since $G$ is a groupoid, by Quillen's
Theorem B. It follows that
$\epsilon$ is a weak equivalence of $G$-diagrams.
\end{proof}

It is shown in \cite[VI.4.2 (p.330)]{GJ} that the homotopy colimit
functor $A \mapsto \hocolim_{G}\ A$ takes pointwise fibrations
to fibrations over $BG$. We know that both the homotopy colimit
functor and the pullback functor $X \mapsto \pb(X)$ preserve weak
equivalences, and so we have the following: 

\begin{lemma}\label{lem 23}
The functors
\begin{equation*}
\hocolim_{G}: s\mathbf{Set}^{G} \leftrightarrows
s\mathbf{Set}/ BG: \pb
\end{equation*}
induce an adjoint equivalence of the associated homotopy categories.
\end{lemma}

Suppose that $\mathcal{M}$ is a right proper closed
model category.  Every morphism $f: X \to Y$ of $\mathcal{M}$ induces
a functor
\begin{equation*}
f^{\ast}: \mathcal{M} / X \to \mathcal{M}
/ Y
\end{equation*}
by composing with $f$. There is a functor
\begin{equation*}
f_{\ast}: \mathcal{M} / Y \to \mathcal{M}
/ X
\end{equation*}
which is defined by pullback along $f$, and $f^{\ast}$ is left adjoint
to $f_{\ast}$. The composition functor plainly preserves cofibrations
and weak equivalences, so the pullback functor preserves fibrations
and trivial fibrations. The pullback functor therefore preserves weak
equivalences of fibrant objects --- note
that a fibrant object of $\mathcal{M} / Y$ is 
a fibration $Z \to Y$.

Each object $\alpha: Z
\to Y$ of $\mathcal{M}/ Y$ has a fibrant model,
meaning a factorization
\begin{equation*}
\xymatrix{
Z \ar[r]^{j_{\alpha}} \ar[dr]_{\alpha} & Z_{\alpha} \ar[d]^{p_{\alpha}} \\
& Y
}
\end{equation*}
where $j_{\alpha}$ is a trivial cofibration and $p_{\alpha}$ is a
fibration. Form the pullback
\begin{equation*}
\xymatrix{
X \times_{Y} Z_{\alpha} \ar[r] \ar[d]_{p_{\alpha\ast}} 
& Z_{\alpha} \ar[d]^{p_{\alpha}} \\
X \ar[r]_{f} & Y
}
\end{equation*}
Then the assignment $\alpha \mapsto p_{\alpha\ast}$ preserves weak
equivalences by the properness assumption for $\mathcal{M}$, and
defines the derived functor
\begin{equation*}
Rf_{\ast}: \Ho(\mathcal{M} / Y) \to \Ho(\mathcal{M}
/ X)
\end{equation*}
Of course, composition with $f$ preserves weak equivalences and
induces a functor
\begin{equation*}
Lf^{\ast}: \Ho(\mathcal{M} / X) \to \Ho(\mathcal{M}
/ Y)
\end{equation*}
Then one shows by chasing explicit homotopy classes that $Lf^{\ast}$ is
left adjoint to $Rf_{\ast}$ . The map $\eta: \beta
\to Rf_{\ast}Lf^{\ast}\beta$ is the map $Z \to X \times_{Y} Z_{f\beta}$
which is determined by the diagram
\begin{equation*}
\xymatrix{
Z \ar[r]^{j_{f\beta}} \ar[d]_{\beta}  & Z_{f\beta} \ar[d]^{p_{f\beta}} \\
X \ar[r]_{f} & Y
}
\end{equation*}
The map $\epsilon: Lf^{\ast}Rf_{\ast}\alpha \to \alpha$ is represented
in the homotopy category, for an object $\alpha: Z \to Y$, by the composite
\begin{equation*}
X \times_{Y} Z_{\alpha} \to Z_{\alpha} \xleftarrow{j_{\alpha}} Z.
\end{equation*}

\begin{lemma}\label{lem 24}
Suppose that $\mathcal{M}$ is a right proper closed
model category, and suppose that $f: X \to Y$ is a weak equivalence of
$\mathcal{M}$. Then the functors
\begin{equation*}
Lf^{\ast}: \Ho(\mathcal{M} / X) \leftrightarrows 
\Ho(\mathcal{M} / Y): Rf_{\ast}
\end{equation*}
form an adjoint equivalence of categories.
\end{lemma}

\begin{proof}
Since $p_{f\beta}$ is a fibration and $f$ is a weak equivalence, the
map $f_{\ast}: X \times_{Y} Z_{f\beta} \to Z_{f\beta}$ is a weak
equivalence. The map $j_{f\beta}: Z \to Z_{f\beta}$ is a weak
equivalence by construction, so that the map $\eta: Z \to X \times_{Y}
Z_{f\beta}$ is a weak equivalence.

Since $p_{\alpha}$ is a fibration and $f$ is a weak equivalence, the
map $f_{\ast}: X \times_{Y} Z_{\alpha} \to Z_{\alpha}$ is a weak
equivalence. It follows that $\epsilon$ is an isomorphism in the
homotopy category.
\end{proof}

The following sequence of results (Corollary \ref{cor 25} -- Corollary
\ref{cor 27}) is perhaps of interest in its own right. It is also a
prototype for a series of results concerning presheaves of groupoids which
appears in the next section.

\begin{corollary}\label{cor 25}
Suppose that the morphism of groupoids $f: G \to H$ induces
a weak equivalence $f: BG \to BH$. Then the composition
with $f$ functor $f^{\ast}$ and the pullback functor $f_{\ast}$
together induce an adjoint equivalence of homotopy categories
\begin{equation*}
Lf^{\ast}: \Ho(s\mathbf{Set}/ BG) \leftrightarrows
\Ho(s\mathbf{Set}/ BH): Rf_{\ast}.
\end{equation*}
\end{corollary}


\begin{corollary}\label{cor 26}
Suppose that the map $f: G \to H$ of groupoids induces a
weak equivalence $f: BG \to BH$. Then the functor 
\begin{equation*}
Rf_{\ast}: \Ho(s\mathbf{Set}^{H}) \to
\Ho(s\mathbf{Set}^{G})
\end{equation*}
which is defined by composition with $f$ is an equivalence of categories. 
\end{corollary}

\begin{proof}
There is a commutative diagram of functors
\begin{equation*}
\xymatrix{
s\mathbf{Set}^{H} \ar[r] \ar[d]_{f_{\ast}}  
& s\mathbf{Set} / BH \ar[d]^{f_{\ast}} \\
s\mathbf{Set}^{G} \ar[r]
& s\mathbf{Set} / BG
}
\end{equation*}
where the horizontal functors are defined by homotopy colimit, and
hence induce equivalences of homotopy categories according to Lemma
\ref{lem 23}.  The functor
\begin{equation*}
f_{\ast}: s\mathbf{Set}/ BH
\to s\mathbf{Set}/ BG
\end{equation*}
is defined by pullback along the map $f: BG \to BH$, and hence
induces an equivalence of homotopy categories by
Corollary \ref{cor 25}
\end{proof}

The restriction functor $f_{\ast}: s\mathbf{Set}^{H} \to
s\mathbf{Set}^{G}$ has a left adjoint $f^{\ast}$ defined by left
Kan extension. The functor $f_{\ast}$ preserves pointwise weak
equivalences and pointwise fibrations, so that the functor $f^{\ast}$
preserves cofibrations and trivial cofibrations, and thus preserves
pointwise weak equivalences between cofibrant objects. It follows that
if $CX$ denotes a cofibrant replacement for a diagram
$X$ on the groupoid $G$, then the assignment $X \to
f^{\ast}CX$ induces a functor
\begin{equation*}
Lf^{\ast} : \Ho(s\mathbf{Set}^{G}) \to \Ho(s\mathbf{Set}^{H})
\end{equation*}
which is left adjoint to the functor
\begin{equation*}
Rf_{\ast}: \Ho(s\mathbf{Set}^{H}) \to \Ho(s\mathbf{Set}^{G})
\end{equation*}
The functor $Rf_{\ast}$ is part of an equivalence on the homotopy
category level, with inverse $G$, say. But every equivalence of
categories is an adjoint equivalence \cite[p.93]{M}, so that
$Lf^{\ast}$ is naturally isomorphic to $G$ as a functor
$\Ho(s\mathbf{Set}^{G}) \to \Ho(s\mathbf{Set}^{H})$.
We have therefore proved the following:

\begin{corollary}\label{cor 27}
Suppose that $f: G \to H$ is a morphism of groupoids such
that $f: BG \to BH$ is a weak equivalence of simplicial
sets. Then the left Kan extension $f^{\ast}$ of the restriction
functor $f_{\ast}: s\mathbf{Set}^{H} \to s\mathbf{Set}^{G}$
has a derived functor
\begin{equation*}
Lf^{\ast} : \Ho(s\mathbf{Set}^{G}) \to \Ho(s\mathbf{Set}^{H})
\end{equation*}
which is an inverse up to natural isomorphism for the derived
restriction functor
\begin{equation*}
Rf_{\ast}: \Ho(s\mathbf{Set}^{H}) \to \Ho(s\mathbf{Set}^{G}).
\end{equation*}
\end{corollary}

Here's a result that is well known \cite{Q}, but stated and proved
in a completely functorial manner. We will need the functoriality for
a corresponding result on presheaves of categories which will be used
in the next section of this paper.

\begin{lemma}\label{lem 28}
There are canonical natural weak equivalences $BC^{op} \simeq dX(C)
\simeq BC$ for a suitably defined natural simplicial set $dX(C)$.
\end{lemma}

\begin{proof}
The simplicial set $BC^{op}$ has $n$-simplices given by strings of
arrows
\begin{equation*}
b_{0} \leftarrow b_{1} \leftarrow \dots \leftarrow b_{n}
\end{equation*}
with simplicial structure maps defined in the obvious way. Consider
the bisimplicial set $X(C)$ with $(m,n)$-bisimplices given by all
strings
\begin{equation*}
b_{m} \to \dots \to b_{0} \to a_{0} \to \dots \to a_{n}.
\end{equation*}
Assigning the $m$-simplex 
\begin{equation*}
b_{m} \to \dots \to b_{0}
\end{equation*}
to this bisimplex defines a function $\phi: X(C)_{m,n} \to
BC^{op}_{m}$, and this list of functions defines a bisimplicial set
map $\phi: X(C) \to BC^{op}$. Assigning the $n$-simplex
\begin{equation*}
a_{0} \to \dots \to a_{n}.
\end{equation*}
to the same bisimplex defines a function $\psi: X(C)_{m,n} \to
BC_{n}$, and the list of functions defines a bisimplicial set map
$\phi: X(C) \to BC$. 

The fibre of $\psi: X(C)_{\ast,n} \to BC_{n}$
over a fixed $n$-simplex $a_{0} \to \dots \to a_{n}$ can be identified
with the simplicial set $B(a_{0} / C^{op})$, which is
contractible. It follows that $\psi$ induces a weak equivalence of
associated diagonal simplicial sets.

The fibre of $\phi: X(C)_{m,\ast} \to BC^{op}_{m}$ over a fixed
$m$-simplex $b_{m} \to \dots \to b_{0}$ can be identified with the
simplicial set $B(b_{0} / C)$, which is again
contractible. It follows that $\phi$ induces a weak equivalence of
associated diagonal simplicial sets.

We have therefore constructed natural weak equivalences
\begin{equation*}
BC^{op} \xleftarrow{\phi} dX(C) \xrightarrow{\psi} BC,
\end{equation*}
as required. Here, $d$ denotes the diagonal functor.
\end{proof}

\section{Presheaves of groupoids}

Suppose that $G$ is a presheaf of groupoids, and let $\mathcal{C}
/ G$ be the corresponding site fibred over $G$. Recall from
Lemma \ref{lem 12} that a
presheaf on $\mathcal{C} / G$ can be identified with an
enriched diagram on the presheaf of groupoids $G^{op}$. It follows
that a simplicial presheaf on $\mathcal{C} / G$ can be
identified with an enriched diagram $X$ on $G^{op}$ taking values in
simplicial sets. 

This means that $X$ consists of functors $X(U):
G(U)^{op} \to s\mathbf{Set}$, $x \mapsto X(U)_{x}$, one for each
object $U \in \mathcal{C}$, such that each morphism $\phi: V \to U$ of
$\mathcal{C}$ induces simplicial set maps $\phi^{\ast}: X(U)_{x} \to
X(V)_{\phi^{\ast}x}$. In addition we require the diagram of simplicial
sets
\begin{equation*}
\xymatrix{
X(U)_{x} \ar[r]^{\alpha_{\ast}} \ar[d]_{\phi^{\ast}} 
& X(U)_{y} \ar[d]^{\phi^{\ast}} \\
X(V)_{\phi^{\ast}(x)} \ar[r]_{(\phi^{\ast}(\alpha))_{\ast}} 
& X(V)_{\phi^{\ast}(y)}
}
\end{equation*} 
to commute for each morphism $\alpha : x \to y$ of $G(U)^{op}$. 

Bundling the simplicial sets $X(U)_{x}$ together over $\Ob(G^{op}(U))$
for all $U$ defines the object map $X_{0} \to \Ob(G^{op})$ of
simplicial presheaves. Recall that $X_{0} \to \Ob(G^{op}) = \Ob(G)$
represents the simplicial presheaf $\psi_{\ast}X$, where $\psi: \Ob(G)
\to G$ is the canonical functor. Lemma \ref{lem 14} implies that a map
$f: X \to Y$ is a local weak equivalence of enriched $G^{op}$-diagrams
if and only if the corresponding map
\begin{equation*}
\xymatrix@C=5pt{
X_{0} \ar[rr]^{f} \ar[dr] && Y_{0} \ar[dl] \\
& \Ob(G^{op})
}
\end{equation*}
is a local weak equivalence of simplicial presheaves on the fibred
site $\mathcal{C} / \Ob(G^{op})$. A similar observation
holds for monomorphisms: a map $g: A \to B$ is a monomorphism of
enriched diagrams if and only if the object-level map $A_{0} \to
B_{0}$ is a monomorphism of simplicial presheaves. Note that Lemma
\ref{lem 15} implies that a global fibration $p: X \to Y$ of enriched
diagrams is an object level global fibration $X_{0} \to Y_{0}$ which
is $G^{op}$-equivariant.

Given an enriched diagram $X$, taking homotopy colimits in each
section defines an enriched homotopy colimit $\hocolim_{G^{op}}\ X$
and a canonical map of simplicial presheaves
\begin{equation*}
\pi: \hocolim_{G^{op}}\ X \to BG^{op}
\end{equation*}
Conversely, one can start with a map $f: Y \to BG^{op}$ and produce an
enriched $G^{op}$-diagram $\pb(Y)$: one applies the
construction which associates the $G^{op}(U)$ diagram $\pb(Y(U))$
to the simplicial set map $Y(U) \to BG^{op}(U)$ in each
section. By working section by section, one sees that there are
natural maps $\eta: Y \to \hocolim_{G^{op}}\ \tilde{Y}$ and $\epsilon:
\pb(\hocolim_{G^{op}}\ X) \to X$, and that these two maps satisfy
the triangle identities. We know from Corollary \ref{cor 21} that the
map $\eta$ is a sectionwise weak equivalence. Corollary \ref{cor 22}
says that all maps 
\begin{equation*}
\epsilon_{x}: \pb(\hocolim_{G^{op}(U)}\ X)_{x}
\to X_{x}
\end{equation*}
are weak equivalences of simplicial sets, for all $x \in
\Ob(G^{op}(U))$ and all $U \in \mathcal{C}$. It follows that $\epsilon$
is a natural weak equivalence of simplicial presheaves on $\mathcal{C}
/ \Ob(G^{op})$.

Lemma \ref{lem 16} says that the homotopy colimit functor
\begin{equation*}
s\Pre(\mathcal{C} / G) \to s\Pre(\mathcal{C}) /
BG^{op} 
\end{equation*}
preserves local weak equivalences. 

\begin{lemma}\label{lem 29}
The functor $s\Pre(\mathcal{C}) / BG^{op} \to
s\Pre(\mathcal{C} / G)$ defined by $X \mapsto \pb(X)$
preserves local weak equivalences.
\end{lemma}

\begin{proof}
Suppose that $Y \to B\Gamma$ is a simplicial set map, where $\Gamma$
is a groupoid. Quillen's Theorem B implies that the square portion of
the diagram
\begin{equation*}
\xymatrix{
\bigsqcup_{x \in \Ob(\Gamma)}\ \pb(Y)_{x} \ar[r] \ar[d] & Y \ar[d] \\
\bigsqcup_{x \in \Ob(\Gamma)}\ B(\Gamma / x) \ar[d] \ar[r] &
B\Gamma \\
\Ob(\Gamma) 
}
\end{equation*}
is homotopy cartesian. Applying this construction in each section to a
simplicial presheaf map $X \to BG^{op}$ gives a diagram of simplicial
presheaf maps
\begin{equation*}
\xymatrix{
\pb(X)_{0} \ar[r] \ar[d] & X \ar[d] \\
\pb(BG^{op})_{0} \ar[d] \ar[r] &
BG^{op} \\
\Ob(G^{op})
}
\end{equation*}
in which the square is homotopy cartesian. Thus if $f: X \to Y$ is a
local weak equivalence of simplicial presheaves over $BG^{op}$, the
induced map $\pb(X)_{0} \to \pb(Y)_{0}$ is a local weak
equivalence of simplicial presheaves over $\Ob(G^{op})$.
The desired statement is then a consequence of Lemma \ref{lem 14}.
\end{proof}

We have assembled a proof of the following:

\begin{theorem}\label{thm 30}
Suppose that $G$ is a presheaf of groupoids on a site $\mathcal{C}$. 
Then the homotopy colimit and pullback functors determine an adjoint
equivalence
\begin{equation*}
\hocolim: \Ho(s\Pre(\mathcal{C} / G)) \simeq \Ho(s\Pre(\mathcal{C})
/ BG^{op}): \pb
\end{equation*}
\end{theorem}

We now have a list of corollaries which is analogous to the sequence
Corollary \ref{cor 25} --- Corollary \ref{cor 27}.

\begin{corollary}\label{cor 31}
Suppose that the map $f: G \to H$ of presheaves of groupoids induces a
local weak equivalence $f: BG \to BH$. Then the derived functor
\begin{equation*}
Rf_{\ast}: \Ho(s\Pre(\mathcal{C}/ H)) \to
\Ho(s\Pre(\mathcal{C} / G))
\end{equation*}
defined by composition with $f$ is an equivalence of categories. 
\end{corollary}

\begin{proof}
There is a commutative diagram of functors
\begin{equation*}
\xymatrix{
s\Pre(\mathcal{C} / H) \ar[r] \ar[d]_{f_{\ast}}  
& s\Pre(\mathcal{C}) / BH^{op} \ar[d]^{f_{\ast}} \\
s\Pre(\mathcal{C} / G) \ar[r]
& s\Pre(\mathcal{C}) / BG^{op}
}
\end{equation*}
where the horizontal functors are defined by homotopy colimit, and
hence induce equivalences of homotopy categories according to Theorem
\ref{thm 30}. The functor 
\begin{equation*}
f_{\ast}: s\Pre(\mathcal{C}) / BH^{op} 
\to s\Pre(\mathcal{C}) / BG^{op}
\end{equation*}
is defined by pullback along the map $f: BG^{op} \to BH^{op}$. This
map $f$ is a local weak equivalence by Lemma \ref{lem 28}, and
pullback along $f: BG^{op} \to BH^{op}$ induces an equivalence of
homotopy categories by Lemma \ref{lem 24}.
\end{proof}

Recall (see the remarks following Lemma \ref{lem 17}) that the
derived functor
\begin{equation*}
Rf_{\ast}: \Ho(s\Pre(\mathcal{C} / H) \to
\Ho(s\Pre{C} / G) 
\end{equation*}
has a left adjoint
\begin{equation*}
Lf^{\ast}: \Ho(s\Pre(\mathcal{C} / G) \to
\Ho(s\Pre{C} / H)
\end{equation*}
which is the homotopy left Kan extension with respect to the
projective model structure.  Corollary \ref{cor 31} implies that the
derived functors $Rf_{\ast}$ and $Lf^{\ast}$ therefore determine an
adjoint equivalence of homotopy categories, and so we have proved the
following:

\begin{corollary}\label{cor 32}
Suppose that $f: G\to H$ is a morphism of presheaves of groupoids such
that $f: BG \to BH$ is a local weak equivalence of simplicial
presheaves. Then the left Kan extension $f^{\ast}$ of the restriction
functor $f_{\ast}: s\Pre(\mathcal{C} / H) \to
s\Pre(\mathcal{C} / G)$ has a derived functor
\begin{equation*}
Lf^{\ast} : \Ho(s\Pre(\mathcal{C} / G) \to
\Ho(s\Pre{C} / H)
\end{equation*}
which is an inverse up to natural isomorphism for the derived
restriction functor
\begin{equation*}
Rf_{\ast}: \Ho(s\Pre(\mathcal{C} / H) \to
\Ho(s\Pre{C} / G)
\end{equation*}
\end{corollary}

Corollary \ref{cor 32} says that the Quillen adjunction determined by
the functor $f: G \to H$ is a Quillen equivalence if $f: BG \to BH$ is
a weak equivalence. The following is essentially a reformulation of
that statement.

\begin{corollary}\label{cor 33}
Suppose that $f: G \to H$ is a morphism of presheaves of groupoids which
induces a local weak equivalence $BG \to BH$. Then the following
statements hold:
\begin{itemize}
\item[1)] 
Suppose that $X$ is a projective cofibrant enriched $G$-diagram and that
$\alpha: f^{\ast}X \to Ff^{\ast}X$ is a weak equivalence of enriched
$H$-diagrams with $Ff^{\ast}X$ projective fibrant. Then the composite
\begin{equation*}
X \xrightarrow{\eta} f_{\ast}f^{\ast}X \xrightarrow{f_{\ast}\alpha}
f_{\ast}Ff^{\ast}X
\end{equation*}
is a weak equivalence of enriched $G$-diagrams.
\item[2)] Suppose that $Y$ is a projective fibrant enriched
$H$-diagram and that $\beta: Cf_{\ast}Y \to f_{\ast}Y$ is a weak
equivalence of enriched $G$-diagrams with $Cf_{\ast}Y$ projective
cofibrant. Then the composite
\begin{equation*}
f^{\ast}Cf_{\ast}Y \xrightarrow{f^{\ast}\beta} f^{\ast}f_{\ast}Y
\xrightarrow{\epsilon} Y
\end{equation*}
is a weak equivalence of enriched $H$-diagrams.
\end{itemize}
\end{corollary}

The following result (Lemma \ref{lem 35}) requires an independent
proof, because the terminal object $\ast$ of $s\Pre(\mathcal{C}
/ G)$ is not projective cofibrant in general.

\begin{example}
Suppose that $K$ is a group. Then $K$ acts freely on the space $EK$
and the map $EK \to \ast$ is a $K$-equivariant trivial fibration,
while a $K$-equivariant map $\ast \to EK$ would pick out a fixed
point. There are no such fixed points, and it follows that the trivial
projective fibration $EK \to \ast$ does not have a section.
\end{example}
 
\begin{lemma}\label{lem 35}
Suppose that $f: G \to H$ is a morphism of presheaves of groupoids
such that the induced map $BG \to BH$ is a weak equivalence of
simplicial presheaves. Then the canonical map $f^{\ast}(\ast) \to
\ast$ is a local weak equivalence of simplicial presheaves on
$\mathcal{C} / H$.
\end{lemma}

\begin{proof}
For a fixed object $U \in \mathcal{C}$, the (simplicial) set
$f^{\ast}(\ast)$ is defined for $y \in H^{op}(U)$ by the assignment
\begin{equation*}
f^{\ast}(\ast)(y) = \varinjlim_{f(x) \to y}\ \ast
\end{equation*}
where the colimit is computed over the index category $f /
x$, and where $f: G(U)^{op} \to H(U)^{op}$ is the corresponding
groupoid morphism. In other words, there is a natural isomorphism
\begin{equation*}
f^{\ast}(y) \cong \pi_{0}B(f / y).
\end{equation*}
Each diagram
\begin{equation*}
\xymatrix{
\bigsqcup_{y \in \Ob(H)^{op}(U)}\ B(f / y) \ar[r] \ar[d]
& BG(U)^{op} \ar[d]^{f} \\
\bigsqcup_{y \in \Ob(H)^{op}(U)}\ B(H^{op}(U) / y) \ar[r] & BH(U)^{op}
}
\end{equation*}
is homotopy cartesian by Quillen's Theorem B, and so the diagram of
simplicial presheaf maps
\begin{equation*}
\xymatrix{
\pb(BG^{op})_{0} \ar[r] \ar[d] & BG^{op} \ar[d]^{f} \\
\pb(BH^{op})_{0} \ar[r] & BH^{op}
}
\end{equation*}
is homotopy cartesian. The simplicial presheaf map $f$ is a weak
equivalence by assumption, and so it follows that there are local weak
equivalences
\begin{equation*}
\pb(BG^{op})_{0} \xrightarrow{\simeq} \pb(BH^{op})_{0}
\xrightarrow{\simeq} \Ob(H^{op}).
\end{equation*}
In particular the presheaf map
\begin{equation*}
\pi_{0}\pb(BG^{op})_{0} \to \Ob(H^{op})
\end{equation*}
induces an isomorphism of associated sheaves. But we also know that
there is an isomorphism
\begin{equation*}
\pi_{0}\pb(BG^{op})_{0} \cong f^{\ast}(\ast)_{0}, 
\end{equation*}
and the resulting
map
\begin{equation*}
f^{\ast}(\ast)_{0} \to \Ob(H^{op})
\end{equation*}
is induced by the canonical morphism $f^{\ast}(\ast) \to \ast$, and it
follows that this canonical morphism is a weak equivalence.
\end{proof}

Write $s_{\ast}\Pre(\mathcal{C} / G)$ for the category of
pointed simplicial presheaves on the site $\mathcal{C} /
G)$. Pointed simplicial presheaves $X$ on $\mathcal{C} / G$
restrict objects $X_{0} \to \Ob(G)$ with a fixed choice of section $s:
\Ob(G) \to X_{0}$ and one can work with this internally, but it's much
easier to work directly with the restriction functor
\begin{equation*}
\psi_{\ast}: s_{\ast}\Pre(\mathcal{C} / G) \to
s_{\ast}\Pre(\mathcal{C} / \Ob(G)).
\end{equation*} 
A similar remark can be made about presheaves of spectra.

The functor $f_{\ast}$ restricts to a functor
\begin{equation*}
f_{\ast}: s_{\ast}\Pre(\mathcal{C}/ H) \to
s_{\ast}\Pre(\mathcal{C} / G)
\end{equation*}
relating pointed simplicial presheaves for the two sites. The functor
$f_{\ast}$ has a left adjoint 
\begin{equation*}
\tilde{f}^{\ast}: s_{\ast}\Pre(\mathcal{C}/ G) \to
s_{\ast}\Pre(\mathcal{C} / H)
\end{equation*}
which is defined for a pointed simplicial presheaf $X$ by 
\begin{equation*}
\tilde{f}^{\ast}(X) = f^{\ast}(X)/f^{\ast}(\ast)
\end{equation*}
 
Fibrant models are formed in pointed simplicial presheaves just as in
simplicial presheaves, and we know from Lemma \ref{lem 35} that the
map $f^{\ast}(\ast) \to \ast$ is a weak equivalence if $f: G \to H$
induces a weak equivalence $BG \to BH$. Suppose that $BG \to BH$ is a
weak equivalence, and suppose that $X$ is a projective cofibrant
pointed simplicial presheaf on $\mathcal{C} / G$. Then, in
the diagram
\begin{equation*}
\xymatrix{
X \ar[r]^{\eta} \ar[dr]_{\eta} 
& f_{\ast}f^{\ast}X \ar[r]^{f_{\ast}(\alpha)} \ar[d] 
& f_{\ast}Ff^{\ast}X \ar[d] \\
& f_{\ast}\tilde{f}^{\ast}X \ar[r]_{f_{\ast}(\alpha)} 
& f_{\ast}F\tilde{f}^{\ast}X
}
\end{equation*}
the map $f_{\ast}Ff^{\ast}X \to f_{\ast}F\tilde{f}^{\ast}X$ is a weak
equivalence, so Corollary \ref{cor 33} implies that the bottom composite
\begin{equation*}
X \xrightarrow{\eta} f_{\ast}\tilde{f}^{\ast}X
\xrightarrow{f_{\ast}(\alpha)} f_{\ast}F\tilde{f}^{\ast}X
\end{equation*}
is a weak equivalence if $\alpha: \tilde{f}^{\ast}X \to
F\tilde{f}^{\ast}X$ is a projective fibrant model for
$\tilde{f}^{\ast}X$.

Suppose that $Y$ is a projective fibrant pointed simplicial presheaf on
$\mathcal{C} / H$. Form the diagram
\begin{equation*}
\xymatrix{
f_{\ast}F\tilde{f}^{\ast}Cf_{\ast}Y \ar[r]^{f_{\ast}F\tilde{f}^{\ast}\beta} 
& f_{\ast}F\tilde{f}^{\ast}f_{\ast}Y \ar[r]^{f_{\ast}F\epsilon}
& f_{\ast}FY \\
f_{\ast}\tilde{f}^{\ast}Cf_{\ast}Y \ar[u]^{f_{\ast}\alpha} 
\ar[r]^{f_{\ast}\tilde{f}^{\ast}\beta}
& f_{\ast}\tilde{f}^{\ast}f_{\ast}Y \ar[u]^{f_{\ast}\alpha} 
\ar[r]^{f_{\ast}\epsilon}
& f_{\ast}Y \ar[u]_{f_{\ast}\alpha} \\
Cf_{\ast}Y \ar[u]^{\eta} \ar[r]_{\beta} & f_{\ast}Y \ar[u]_{\eta} \ar[ur]_{1}
}
\end{equation*}
by making suitable choices of fibrant models $\alpha$ and cofibrant
models $\beta$. Then the composite
\begin{equation*}
Cf_{\ast}Y \xrightarrow{\eta} f_{\ast}\tilde{f}^{\ast}Cf_{\ast}Y
\xrightarrow{f_{\ast}\alpha} f_{\ast}F\tilde{f}^{\ast}Cf_{\ast}Y
\end{equation*}
is a weak equivalence from what we have just seen, since $Cf_{\ast}Y$
is projective cofibrant. The map $f_{\ast}\alpha: f_{\ast}Y \to
f_{\ast}FY$ is a weak equivalence by Corollary \ref{cor 18}, and of
course the map $\beta: Cf_{\ast}Y \to f_{\ast}Y$ is a weak
equivalence. It follows that the top composite in the diagram is a
weak equivalence. The composite map
\begin{equation*}
\tilde{f}^{\ast}Cf_{\ast}Y \xrightarrow{\tilde{f}^{\ast}\beta}
\tilde{f}^{\ast}f_{\ast}Y \xrightarrow{\epsilon} Y
\end{equation*}
is therefore a weak equivalence on account of Corollary \ref{cor 32},
since $f_{\ast}$ must then reflect weak equivalences between
projective fibrant objects.

We have therefore proved the following:

\begin{lemma}\label{lem 43}
Suppose that $f: G \to H$ is a morphism of presheaves of groupoids which
induces a local weak equivalence $BG \to BH$. Then the following
statements hold:
\begin{itemize}
\item[1)] Suppose that $X$ is a projective cofibrant pointed
simplicial presheaf on $\mathcal{C} / G$ and that $\alpha:
f^{\ast}X \to Ff^{\ast}X$ is a weak equivalence of pointed simplicial
presheaves on $\mathcal{C} / H$ with $Ff^{\ast}X$ projective
fibrant. Then the composite
\begin{equation*}
X \xrightarrow{\eta} f_{\ast}\tilde{f}^{\ast}X \xrightarrow{f_{\ast}\alpha}
f_{\ast}F\tilde{f}^{\ast}X
\end{equation*}
is a weak equivalence of pointed simplicial presheaves on $\mathcal{C}
/ G$.
\item[2)] Suppose that $Y$ is a projective fibrant pointed simplicial
presheaf on $\mathcal{C} / H$ and that $\beta:
C\tilde{f}_{\ast}Y \to \tilde{f}_{\ast}Y$ is a weak equivalence of
pointed simplicial presheaves on $\mathcal{C} / G$ with
$C\tilde{f}_{\ast}Y$ projective cofibrant. Then the composite
\begin{equation*}
\tilde{f}^{\ast}Cf_{\ast}Y \xrightarrow{\tilde{f}^{\ast}\beta} 
\tilde{f}^{\ast}f_{\ast}Y
\xrightarrow{\epsilon} Y
\end{equation*}
is a weak equivalence of pointed simplicial presheaves on $\mathcal{C}
/ H$.
\end{itemize}
\end{lemma}

\begin{corollary}\label{cor 37}
Suppose that $f: G\to H$ is a morphism of presheaves of groupoids such
that $f: BG \to BH$ is a local weak equivalence of simplicial
presheaves. Then the left Kan extension $\tilde{f}^{\ast}$ of the restriction
functor $f_{\ast}: s_{\ast}\Pre(\mathcal{C} / H) \to
s_{\ast}\Pre(\mathcal{C} / G)$ has a derived functor
\begin{equation*}
L\tilde{f}^{\ast} : \Ho(s_{\ast}\Pre(\mathcal{C} / G)) \to
\Ho(s_{\ast}\Pre{C} / H))
\end{equation*}
which is an inverse up to natural isomorphism for the derived
restriction functor
\begin{equation*}
Rf_{\ast}: \Ho(s_{\ast}\Pre(\mathcal{C} / H)) \to
\Ho(s_{\ast}\Pre{C} / G))
\end{equation*}
\end{corollary}

There is one final thing to know about pointed simplicial presheaves
on $\mathcal{C} / G$, which will be of some use later on:

\begin{lemma}\label{lem 45}
Suppose that $K$ is a pointed simplicial set. Then the following hold:
\begin{itemize}
\item[1)] 
If $p: X \to Y$ is a projective fibration (respectively
trivial projective fibration) then the induced map of pointed function
complexes 
\begin{equation*}
p_{\ast}: \mathbf{hom}_{\ast}(K,X) \to
\mathbf{hom}_{\ast}(K,Y) 
\end{equation*}
is a projective fibration (respectively
trivial projective fibration).
\item[2)]
The functor $X \mapsto X \wedge K$ preserves projective cofibrations 
and trivial projective cofibrations.
\end{itemize}
\end{lemma}

\begin{proof}
Statement 1) follows from the fact that restriction along the functor
$\psi:\Ob(G) \to G$ preserves the displayed pointed function complex
constructions. Statement 2) is equivalent to Statement 1), by an
adjointness argument.
\end{proof}

Suppose again that $G$ is a presheaf of groupoids on the site
$\mathcal{C}$, and write $\mathbf{Spt}(\mathcal{C} / G)$ for
the category of presheaves of spectra on the fibred
site $\mathcal{C} / G$. 

Let $\psi: \Ob(G) \to G$ denote the canonical functor, and recall
that a map $f: X \to Y$ of pointed simplicial presheaves is a local
weak equivalence (respectively cofibration) if and only if its
restriction $f_{\ast}: \psi_{\ast}X \to \psi_{\ast}Y$ is a local weak
equivalence (respectively cofibration) on the site $\mathcal{C}
/ \Ob(G)$. By definition, $f$ is a projective fibration if
and only if $f_{\ast}$ is a global fibration on $\mathcal{C}
/ \Ob(G)$. We also know, from Lemma \ref{lem 15}, that the
restriction functor $\psi_{\ast}$ preserves global fibrations.

Recall \cite{J3} that a map $g: Z \to W$ of presheaves of spectra is
a stable equivalence if the induced map $QJX \to QJY$ is a levelwise
weak equivalence, where $X \to JX$ is a natural choice of strictly
fibrant model and $QY$ for a level fibrant object $Y$ is the result of
the usual stabilization construction. In particular, $QY^{n}$ is the
colimit of the diagram
\begin{equation*}
Y^{n} \to \Omega Y^{n+1} \to \Omega^{2}Y^{n+2} \to \cdots
\end{equation*}
Restriction along the canonical functor $\psi: \Ob(G) \to G$ preserves
level fibrant models and the stabilization construction (the latter by Lemma
\ref{lem 35}). The restriction functor also reflects level weak
equivalences. It follows that a map $g: Z \to W$ of presheaves of
spectra on the site $\mathcal{C} / G$ is a stable equivalence
if and only if its restriction $g_{\ast}: \psi_{\ast}Z \to
\psi_{\ast}W$ is a stable equivalence of presheaves of spectra on the
site $\mathcal{C} / \Ob(G)$. It is also relatively easy to
see that $g$ is a cofibration of presheaves of spectra on $\mathcal{C}
/ A$ if and only if $g_{\ast}$ is a cofibration of presheaves
of spectra on $\mathcal{C} / \Ob(G)$.

It can be shown that a map $p: X \to Y$ of presheaves of spectra is a
stable fibration if and only if the following conditions hold:
\begin{itemize}
\item[1)]
All level maps $p: X^{n} \to Y^{n}$ are fibrations of pointed 
simplicial presheaves.
\item[2)]
Given any commutative diagram
\begin{equation*}
\xymatrix{
X \ar[r]^{j} \ar[d]_{p} & Z \ar[d] \\
Y \ar[r]_{j} & W
}
\end{equation*}
where the maps labelled $j$ are stable equivalences and $Z$ and $W$ 
are stably fibrant, then all diagrams
\begin{equation*}
\xymatrix{
X^{n} \ar[r]^{j} \ar[d]_{p} & Z^{n} \ar[d] \\
Y^{n} \ar[r]_{j} & W^{n}
}
\end{equation*}
are homotopy cartesian diagrams of pointed simplicial presheaves.
\end{itemize}
In particular, if $X$ and $Y$ are already stably fibrant, then a
stable fibration $p: X \to Y$ is a map such that all level maps $p:
X^{n} \to Y^{n}$ are fibrations of pointed simplicial presheaves.

Say that a map $p: X \to Y$ of presheaves of spectra on $\mathcal{C}
/ G$ is a projective fibration if the restriction $p_{\ast}:
\psi_{\ast}X \to \psi_{\ast}Y$ is a stable fibration on $\mathcal{C}
/ \Ob(G)$. One can see by using the criteria 1) and 2) above
that the functor $\psi_{\ast}$ preserves stable fibrations, so that
every stable fibration is a projective fibration.  A projective
cofibration of presheaves of spectra on $\mathcal{C} / G$ is
a map which has the left lifting property with respect to all maps
which are simultaneously stable equivalences and projective
fibrations.

The restriction functor $\psi_{\ast}$ preserves stable fibrations and
trivial stable fibrations, so that its adjoint $\psi^{\ast}$ preserves
cofibrations and stably trivial cofibrations.  The stable model
structure of presheaves of spectra is cofibrantly generated, so we are
therefore entitled to the following analogue of Lemma \ref{lem 17}:.

\begin{lemma}
The category $\mathbf{Spt}(\mathcal{C} / A)$ of presheaves of
spectra on the site $\mathcal{C} / G$, together with the
stable weak equivalences, projective fibrations and projective
cofibrations as defined above, satisfies the axioms for a closed model
category.
\end{lemma} 

A presheaf of spectra $X$ is stably fibrant if all objects
$X^{n}$ are fibrant and all adjoint bonding maps $X^{n} \to \Omega
X^{n+1}$ weak equivalences. Furthermore, a map $f:X \to Y$ between
stably fibrant presheaves of spectra is a stable equivalence if and
only if all level maps $X^{n} \to Y^{n}$ are weak equivalences of
simplicial presheaves. It follows that a presheaf of spectra $Z$ on
$\mathcal{C} / G$ is projective fibrant if and only if all
objects $Z^{n}$ are projective fibrant pointed simplicial presheaves
and all morphisms $Z^{n} \to \Omega Z^{n+1}$ are local weak
equivalences. It also follows that a map $g: Z \to W$ of projective
fibrant presheaves of spectra is a stable weak equivalence if and only
if the restriction $g: \psi_{\ast}Z \to \psi_{\ast}W$ is a stable weak
equivalence of presheaves of spectra on $\mathcal{C} /
\Ob(G)$. It is a further consequence that the restriction functor
\begin{equation*}
f_{\ast}: \mathbf{Spt}(\mathcal{C} / H) \to
\mathbf{Spt}(\mathcal{C} / G)
\end{equation*}
preserves projective fibrant presheaves of spectra, and preserves
stable weak equivalences between projective fibrant presheaves of
spectra.

This characterization gives rise to an obvious recognition principle
for projective fibrations of presheaves of spectra on $\mathcal{C}
/ G$, and implies that a map $q: Z \to W$ between projective
fibrant presheaves of spectra is a projective fibration if and only if
all level maps $p: Z^{n} \to W^{n}$ are projective fibrations. It
follows in particular that for any morphism $f: G \to H$ of presheaves of
groupoids the restriction functor
\begin{equation*}
f_{\ast}: \mathbf{Spt}(\mathcal{C} / H) \to
\mathbf{Spt}(\mathcal{C} / G)
\end{equation*}
preserves projective fibrations and stable equivalences between
projective fibrant objects.  It therefore also follows that the left
adjoint
\begin{equation*}
\tilde{f}_{\ast}: \mathbf{Spt}(\mathcal{C} / G) \to
\mathbf{Spt}(\mathcal{C} / H)
\end{equation*}
preserves projective cofibrations and stable equivalences between
projective cofibrant objects. We are therefore entitled to derived functors
\begin{equation*}
L\tilde{f}^{\ast}: \Ho(\mathbf{Spt}(\mathcal{C}/ G)) \leftrightarrows
\Ho(\mathbf{Spt}(\mathcal{C} / H)): Rf_{\ast}
\end{equation*}
relating the associated stable categories. Furthermore,
$L\tilde{f}^{\ast}$ is left adjoint to $Rf_{\ast}$, with the usual
description of unit and counit.

\begin{lemma}\label{lem 40}
Suppose that $f: G \to H$ is a morphism of presheaves of groupoids which
induces a local weak equivalence $BG \to BH$. Then the following
statements hold:
\begin{itemize}
\item[1)] 
Suppose that $X$ is a projective cofibrant presheaf of spectra on
$\mathcal{C} / G$ and that $\alpha: f^{\ast}X \to Ff^{\ast}X$
is a stable equivalence of presheaves of spectra on $\mathcal{C}
/ H$ with $Ff^{\ast}X$ projective fibrant. Then the composite
\begin{equation*}
X \xrightarrow{\eta} f_{\ast}\tilde{f}^{\ast}X \xrightarrow{f_{\ast}\alpha}
f_{\ast}F\tilde{f}^{\ast}X
\end{equation*}
is a stable weak equivalence.
\item[2)] 
Suppose that $Y$ is a projective fibrant presheaf of spectra
on $\mathcal{C} / H$ that $\beta: C\tilde{f}_{\ast}Y \to
\tilde{f}_{\ast}Y$ is a stable equivalence of presheaves of spectra on
$\mathcal{C} / G$ with $C\tilde{f}_{\ast}Y$ projective
cofibrant. Then the composite
\begin{equation*}
\tilde{f}^{\ast}Cf_{\ast}Y \xrightarrow{\tilde{f}^{\ast}\beta} 
\tilde{f}^{\ast}f_{\ast}Y
\xrightarrow{\epsilon} Y
\end{equation*}
is a stable equivalence.
\end{itemize}
\end{lemma}

\begin{proof}
Suppose that $X$ is a projective cofibrant presheaf of spectra. Then
there is a level equivalence $\pi: \tilde{X} \to X$ where the pointed
simplicial presheaf $\tilde{X}^{0}$ is projective cofibrant and all bonding
maps $S^{1} \wedge \tilde{X}^{n} \to \tilde{X}^{n+1}$ are projective cofibrations. In
particular all pointed simplicial presheaves $\tilde{X}^{n}$ are projective
cofibrant. The construction of $\pi$ is the standard cofibrant
replacement trick, which takes advantage of the fact that a map $p: X
\to Y$ is a projective fibration and a stable equivalence if and only
if all level maps $p: X^{n} \to Y^{n}$ are trivial projective fibrations of
pointed simplicial presheaves. In the diagram
\begin{equation*}
\xymatrix{
\tilde{X} \ar[r]^-{\eta} \ar[d]_{\pi} 
& f_{\ast}\tilde{f}^{\ast}\tilde{X} \ar[r]^{f_{\ast}\alpha} 
& f_{\ast}F\tilde{f}^{\ast}\tilde{X} \ar[d]^{f_{\ast}F\tilde{f}^{\ast}\pi} \\
X \ar[r]_-{\eta} & f_{\ast}\tilde{f}^{\ast}X \ar[r]_{f_{\ast}\alpha} 
& f_{\ast}F\tilde{f}^{\ast}X
}
\end{equation*}
the map $f_{\ast}F\tilde{f}^{\ast}\pi$ is a stable equivalence, since
$\tilde{f}^{\ast}$ preserves stable equivalences between projective
cofibrant objects. The top horizontal composite is a level weak
equivalence by Lemma \ref{lem 43}, and so the bottom horizontal
composite is also a stable equivalence.

Assertion 2) has a similar proof: the cofibrant model $Cf_{\ast}Y$ can
be chosen so that it consists of projective cofibrant pointed
simplicial presheaves in all levels.
\end{proof}

We have also proved the following

\begin{theorem}\label{thm 41}
Suppose that $f: G\to H$ is a morphism of presheaves of groupoids such
that $f: BG \to BH$ is a local weak equivalence of simplicial
presheaves. Then the left Kan extension $\tilde{f}^{\ast}$ of the restriction
functor $f_{\ast}: \mathbf{Spt}(\mathcal{C} / H) \to
\mathbf{Spt}(\mathcal{C} / G)$ has a derived functor
\begin{equation*}
L\tilde{f}^{\ast} : \Ho(\mathbf{Spt}(\mathcal{C} / G)) \to
\Ho(\mathbf{Spt}(\mathcal{C} / H))
\end{equation*}
which is an inverse up to natural isomorphism for the derived
restriction functor
\begin{equation*}
Rf_{\ast}: \Ho(\mathbf{Spt}\Pre(\mathcal{C} / H)) \to
\Ho(\mathbf{Spt}(\mathcal{C} / G)).
\end{equation*}
\end{theorem}

Theorem \ref{thm 41} implies the corresponding result for presheaves
of symmetric spectra rather easily, subject to having appropriate
projective model structures in place.

For a fixed
presheaf of groupoids $G$  the restriction functor
\begin{equation*}
\psi_{\ast}: \mathbf{Spt}_{\Sigma}(\mathcal{C} / G) 
\to \mathbf{Spt}_{\Sigma}(\mathcal{C} / \Ob(G))
\end{equation*} 
between the respective categories of presheaves of symmetric spectra
preserves stable fibrations and trivial stable fibrations (see
\cite{J4}). It follows that its left adjoint $\psi^{\ast}$ preserves
cofibrations and trivial cofibrations. Say that a map $p: X \to Y$ of
symmetric spectra on $\mathcal{C} / G$ is a projective
fibration if the induced map $p_{\ast}: \psi_{\ast}X \to \psi_{\ast}Y$
is a stable fibration of presheaves of symmetric spectra on
$\mathcal{C} / \Ob(G)$. A projective cofibration is a map of
$\mathbf{Spt}_{\Sigma}(\mathcal{C} / G)$ which has the left
lifting property with respect to all maps which are both stable
equivalences and projective fibrations. Note that every stable
fibration of $\mathbf{Spt}_{\Sigma}(\mathcal{C} / G)$ is a
projective fibration, so that every projective cofibration is a
cofibration.

The category of presheaves of
symmetric spectra is cofibrantly generated, and one can prove the
following:

\begin{lemma}
The category $\mathbf{Spt}_{\Sigma}(\mathcal{C} / A)$ of
presheaves of spectra on the site $\mathcal{C} / G$.
together with the stable weak equivalences, projective fibrations and
projective cofibrations as defined above, satisfies the axioms for a
closed model category.
\end{lemma} 

Suppose that $f: G \to H$ is a morphism of presheaves of
groupoids. Then the restriction functor
\begin{equation*}
f_{\ast}: \mathbf{Spt}_{\Sigma}(\mathcal{C} / H) \to
\mathbf{Spt}_{\Sigma}(\mathcal{C} / G)
\end{equation*}
preserves projective fibrations and trivial projective fibrations. It
follows that the left adjoint functor
\begin{equation*}
\tilde{f}^{\ast}: \mathbf{Spt}_{\Sigma}(\mathcal{C} / G) \to
\mathbf{Spt}_{\Sigma}(\mathcal{C} / H)
\end{equation*}
preserves projective cofibrations and trivial projective
cofibrations. As in all other cases, one shows that the
corresponding adjunction
\begin{equation*}
L\tilde{f}^{\ast}: \Ho(\mathbf{Spt}_{\Sigma}(\mathcal{C}/ G)) 
\leftrightarrows
\Ho(\mathbf{Spt}_{\Sigma}(\mathcal{C} / H)): Rf_{\ast}
\end{equation*}
is an adjoint equivalence.

Recall that the forgetful functor $U:
\mathbf{Spt}_{\Sigma}(\mathcal{D}) \to \mathbf{Spt}(\mathcal{D})$ has
a left adjoint $V$, and that these functors form a Quillen equivalence,
for any small site $\mathcal{D}$. In particular, $V$ preserves
cofibrations and trivial cofibrations while $U$ preserves stable
fibrations and trivial stable fibrations, and the corresponding
derived functors
\begin{equation*}
LV: \Ho(\mathbf{Spt}(\mathcal{D})) \leftrightarrows
\Ho(\mathbf{Spt}_{\Sigma}(\mathcal{D})): RU
\end{equation*}
form an adjoint equivalence of categories. In particular, there are
(composite) stable equivalences
\begin{equation*}
VCUY \xrightarrow{V\beta} VUY \xrightarrow{\epsilon} Y
\end{equation*}
for all stably fibrant $Y$ and
\begin{equation*}
X \xrightarrow{\eta} UVX \xrightarrow{U\alpha} UFVX
\end{equation*}
for all cofibrant $X$, where $\beta$ and $\alpha$ are cofibrant and
fibrant models, respectively.

\begin{lemma}
Suppose that the map $f: G \to H$ of presheaves of groupoids induces a
local weak equivalence $BG \to BH$. Then the derived functor
\begin{equation*}
Rf_{\ast}: \Ho(\mathbf{Spt}_{\Sigma}(\mathcal{C} / H))
\to \Ho(\mathbf{Spt}_{\Sigma}(\mathcal{C} / G))
\end{equation*}
is an equivalence of categories.
\end{lemma}

\begin{proof}
The diagram of functors
\begin{equation*}
\xymatrix{
\mathbf{Spt}_{\Sigma}(\mathcal{C} / H) \ar[r]^{U} \ar[d]_{f_{\ast}} 
& \mathbf{Spt}(\mathcal{C} / H) \ar[d]^{f_{\ast}} \\
\mathbf{Spt}_{\Sigma}(\mathcal{C} / G) \ar[r]_{U}
& \mathbf{Spt}(\mathcal{C} / G)
}
\end{equation*}
induces a commutative diagram of right derived functors
\begin{equation*}
\xymatrix{
\Ho(\mathbf{Spt}_{\Sigma}(\mathcal{C} / H)) 
\ar[r]^{RU}_{\simeq} \ar[d]_{Rf_{\ast}} 
& \Ho(\mathbf{Spt}(\mathcal{C} / H)) \ar[d]^{Rf_{\ast}}_{\simeq} \\
\Ho(\mathbf{Spt}_{\Sigma}(\mathcal{C} / G)) \ar[r]_{RU}^{\simeq}
& \Ho(\mathbf{Spt}(\mathcal{C} / G))
}
\end{equation*}
by standard results about symmetric spectra and Theorem \ref{thm 41}.
\end{proof}

The left adjoint $L\tilde{f}^{\ast}$ of $Rf_{\ast}$ must coincide with
the inverse of $Rf_{\ast}$ up to natural isomorphism, and so we have
the following:

\begin{corollary}\label{cor 45}
Suppose that $f: G\to H$ is a morphism of presheaves of groupoids such
that $f: BG \to BH$ is a local weak equivalence of simplicial
presheaves. Then the left Kan extension $\tilde{f}^{\ast}$ of the restriction
functor $f_{\ast}: \mathbf{Spt}_{\Sigma}(\mathcal{C} / H) \to
\mathbf{Spt}_{\Sigma}(\mathcal{C} / G)$ has a derived functor
\begin{equation*}
L\tilde{f}^{\ast} : \Ho(\mathbf{Spt}_{\Sigma}(\mathcal{C} / G)) \to
\Ho(\mathbf{Spt}_{\Sigma}(\mathcal{C} / H))
\end{equation*}
which is an inverse up to natural isomorphism for the derived
restriction functor
\begin{equation*}
Rf_{\ast}: \Ho(\mathbf{Spt}_{\Sigma}\Pre(\mathcal{C} / H)) \to
\Ho(\mathbf{Spt}_{\Sigma}(\mathcal{C} / G)).
\end{equation*}
\end{corollary}

\begin{corollary}
Suppose that $f: G \to H$ is a morphism of presheaves of groupoids which
induces a local weak equivalence $BG \to BH$. Then the following
statements hold:
\begin{itemize}
\item[1)] 
Suppose that $X$ is a projective cofibrant presheaf of symmetric spectra on
$\mathcal{C} / G$ and that $\alpha: f^{\ast}X \to Ff^{\ast}X$
is a stable equivalence of presheaves of symmetric spectra on $\mathcal{C}
/ H$ with $Ff^{\ast}X$ projective fibrant. Then the composite
\begin{equation*}
X \xrightarrow{\eta} f_{\ast}\tilde{f}^{\ast}X \xrightarrow{f_{\ast}\alpha}
f_{\ast}F\tilde{f}^{\ast}X
\end{equation*}
is a stable weak equivalence.
\item[2)] Suppose that $Y$ is a projective fibrant presheaf of
symmetric spectra on $\mathcal{C} / H$ that $\beta:
C\tilde{f}_{\ast}Y \to \tilde{f}_{\ast}Y$ is a stable equivalence of
presheaves of symmetric spectra on $\mathcal{C} / G$ with
$C\tilde{f}_{\ast}Y$ projective cofibrant. Then the composite
\begin{equation*}
\tilde{f}^{\ast}Cf_{\ast}Y \xrightarrow{\tilde{f}^{\ast}\beta} 
\tilde{f}^{\ast}f_{\ast}Y
\xrightarrow{\epsilon} Y
\end{equation*}
is a stable equivalence.
\end{itemize}
\end{corollary}


\catcode`\@=11 \renewenvironment{thebibliography}[1]{
\@xp\section\@xp*\@xp{\refname}%
\normalfont\footnotesize\labelsep
.5em\relax\renewcommand\theenumiv{\arabic{enumiv}}\let\p@enumiv\@empty
\list{\@biblabel{\theenumiv}}{\settowidth\labelwidth{\@biblabel{#1}}%
\leftmargin\labelwidth \advance\leftmargin\labelsep
\usecounter{enumiv}}%
\sloppy \clubpenalty\@M
\widowpenalty\clubpenalty \sfcode`\.=\@m }

\def\@biblabel#1{\@ifnotempty{#1}{[#1]}}
\catcode`\@=\active


\bibliographystyle{amsalpha}

\enddocument